\documentclass[12pt]{article}
\usepackage{amsmath, graphicx, amsfonts,amssymb, calrsfs}
\usepackage{color}

\def\bbR{\mathbb{R}}
\def\r{\rho}
\def\s{\sigma}

\def\e{\epsilon}
\def\cK{\mathcal{K}}
\def\cV{\mathcal{V}}
\def\cM{\mathcal{M}}
\def\L{\Lambda}
\def\E{\mathcal{E}}
\def\cF{\mathcal{F}}
\def\o{\omega}
\def\ball{B^n_2}
\def\polar{K^\circ}
\def\l{\lambda}
\def\cS{\mathcal{S}}

\def\be{\begin{equation}}
\def\ee{\end{equation}}
\def\bea{\begin{eqnarray}}
\def\eea{\end{eqnarray}}
\def\bt{\begin{theorem}}
\def\et{\end{theorem}}
\def\bl{\begin{lemma}}
\def\el{\end{lemma}}
\def\br{\begin{remark}}
\def\er{\end{remark}}
\def\bc{\begin{corollary}}
\def\ec{\end{corollary}}
\def\bd{\begin{definition}}
\def\ed{\end{definition}}
\def\bp{\begin{proposition}}
\def\ep{\end{proposition}}
\newtheorem{theorem}{Theorem}[section]
\newtheorem{lemma}{Lemma}[section]
\newtheorem{remark}{Remark}[section]
\newtheorem{proposition}{Proposition}[section]
\newtheorem{corollary}{Corollary}[section]

\newtheorem{definition}{Definition}[section]
\begin{document}
\title{$L_p$ geominimal surface areas and their inequalities
}

\author{Deping Ye }
\date{}
\maketitle
\begin{abstract}
In this paper, we introduce the $L_p$ geominimal surface area for all $-n\neq p<1$, which extends the classical geominimal surface area ($p=1$) by Petty and the $L_p$ geominimal surface area by Lutwak ($p>1$). Our extension of the $L_p$ geominimal surface area is motivated by recent work on the extension of the $L_p$ affine surface area -- a fundamental object in (affine) convex geometry. We prove some properties for the $L_p$ geominimal surface area and its related inequalities, such as, the affine isoperimetric inequality and a Santal\'{o} style inequality. Cyclic inequalities are established to obtain the monotonicity of the $L_p$ geominimal surface areas. Comparison between the $L_p$ geominimal surface area and the $p$-surface area is also provided.
\end{abstract}
   \section{Introduction and Overview of Results} 
The classical isoperimetric problem asks: what is the minimal area among all convex bodies (i.e., convex compact subsets with nonempty interior) $K\subset \bbR^n$  with volume $1$? The solution of this old problem is now known as the (classical) isoperimetric inequality, namely, {the minimal area is attained at and only at Euclidean balls with volume $1$.}
The isoperimetric inequality is an extremely powerful tool in geometry and
related areas. Note that the classical isoperimetric inequality does not have the ``affine invariant" flavor, because the area may change under linear transformations (even) with unit (absolute value of) determinant. 

\vskip 2mm 
However, many objects in (affine) convex geometry are invariant under invertible linear transformations. A typical example is the Mahler volume product $M(K)=|K||\polar|$, the product of the volume of $K$ and its polar body $\polar$. That is, $M(K)=M(TK)$ for all invertible linear transformations $T$ on $\bbR^n$. One can ask a question for $M(K)$  similar to the classical isoperimetric problem: what is the maximum of $M(K)$ among all convex bodies $K$? The celebrated Blaschke-Santal\'{o} inequality states that $M(K)$ attains its maximum at and only at origin-symmetric ellipsoids if $K$ is assumed to have its centroid at the origin. Such an inequality is arguably more important than the classical isoperimetric inequality. Consequently, the Blaschke-Santal\'{o} inequality has numerous applications in (affine) convex geometry and many other related fields, for instance, in quantum information theory \cite{ASY2013, Szarek2011}. The family of affine isoperimetric inequalities continues to grow (see, for example  \cite{CG, HaberlFranz2009, HaberlFranz20091,  LYZ2, LYZ1}). It is well known that the Blaschke-Santal\'{o} inequality was first discovered and proved by using the classical affine isoperimetric inequality \cite{Santalo1949}, which bounds the classical affine surface area from above in terms of volume.

\vskip 2mm The study of affine surface areas goes back to Blaschke 
in \cite{Bl1} about one hundred years ago, and its $L_p$ counterpart was first introduced by Lutwak in \cite{Lu1}. The notion of $L_p$ affine surface area was further
extended to all $p\in \bbR$ for general convex bodies in, e.g.,
\cite{MW2, SW5, SW4}. In fact, extensions of $L_p$ affine surface area to all $-n\neq p\in \bbR$ were obtained by their integral expressions (see e.g.\ Theorem \ref{equivalent:affine:surface:area}) and by 
investigating the asymptotic behavior of the volume of certain
families of convex bodies \cite{Lei1986, MW1, MW2, SW5, SW4, SW1, W1, WY2008} (and even star-shape bodies \cite{WernerYe2010}). 
The $L_p$ affine surface area is now thought to be at the core of
the rapidly developing $L_p$-Brunn-Minkowski theory. It has many nice properties, such as, affine invariance (i.e., invariance under all linear transformations $T$ with $|det(T)|=1$), the valuation property,  semi-continuity (upper for $p>0$ and lower for $-n\neq p<0$), etc.  Moreover, the $L_p$ affine surface area for $p>0$ has been proved to be, roughly speaking, the unique valuation with properties of affine invariance and upper semi-continuity \cite{LR1, LudR}. Recently, it has been connected with the information theory of convex bodies (see e.g.,
\cite{Jenkinson2012, Paouris2010, Werner2012a, Werner2012b}). Affine
isoperimetric inequalities related to the $L_p$ affine surface area can be found in, for instance, \cite{Lu1,
WY2008}.

\vskip 2mm Another fundamental concept in convex geometry is the (classical) geominimal surface area introduced by Petty \cite{Petty1974} about 40 years ago. As explained by Petty in \cite{Petty1974}, the (classical) geominimal surface area naturally connects affine geometry, relative geometry and Minkowski geometry. Hence it received a lot of attention (e.g., \cite{Petty1974, Petty1985, Schneider2013}). In fact, it can be viewed as a sibling concept of the classical affine surface area as their definitions are similar. However, they are in general different from each other; for example the (classical) geominimal surface area is continuous while the classical affine surface area is only upper semi-continuous on the set of all convex bodies (equipped with the Hausdorff metric). 

\vskip 2mm  In his seminal paper \cite{Lu1}, Lutwak introduced $L_p$  geominimal surface area for $p> 1$.  $L_p$ geominimal surface area has many properties very similar to $L_p$ affine surface area, for instance, affine invariance with the same degree of homogeneous. For certain class of convex bodies, the $L_p$ geominimal surface area for $p\geq 1$ is equal to the $L_p$ affine surface area (see \cite{Lu1} or Proposition \ref{geominimal:ball-1}). However, the $L_p$ geominimal surface area is different from the $L_p$ affine surface area, because the former one is continuous while the latter one is only upper semi-continuous on the set of all convex bodies (equipped with the Hausdorff metric). Moreover, as mentioned before, there are nice integral expressions for the $L_p$ affine surface area (see e.g. Theorem \ref{equivalent:affine:surface:area}), while finding similar integral expressions for the $L_p$ geominimal surface area appears to be impossible. Note that the extension of $L_p$ affine surface area from $p\geq 1$ to all $-n\neq p\in \bbR$ mainly relied on the integral expression of $L_p$ affine surface area. Recently, the $L_p$ affine surface area was further extended to its Orlicz counterpart, general affine surface areas, involving general (even nonhomogeneous) convex or concave functions \cite{Ludwig2009, LR1}.

\vskip 2mm  Affine isoperimetric inequalities related to the $L_p$ geominimal surface area were proved for $p=1$ in \cite{Petty1974, Petty1985} and for  $p> 1$ in \cite{Lu1}. Roughly speaking, those affine isoperimetric inequalities assert that among all convex bodies with fixed volume and with centroid at the origin (the condition on the centroid can be removed for $p=1$), the $L_p$ geominimal surface area attains its maximum at and only at origin-symmetric ellipsoids. As mentioned in \cite{Lu1}, the affine isoperimetric inequality related to the $L_p$ geominimal surface area for $p\geq 1$ is equivalent to the Blaschke-Santal\'{o} inequality in the following sense: either of them can be easily obtained from the other. A Santal\'{o} style inequality for the $L_p$ geominimal surface area was proved in \cite{Zhou2011}.

\vskip 2mm This paper is dedicated to further extend Lutwak's $L_p$ geominimal surface area from $p\geq 1$ to all $-n\neq p\in \bbR$. Hereafter, we use $\tilde{G}_p(K)$ to denote the $L_p$ geominimal surface area of a convex body $K$ (with the origin in its interior). The following Santal\'{o} style inequality will be proved. 

\vskip 2mm \noindent {\bf Theorem \ref{geominimal product:bound} }  {\em Let $K$ be a convex body with centroid (or the Santal\'{o} point) at the origin.  \\ (i). For $p\geq 0$,  $$\tilde{G}_p(K)\tilde{G}_p(\polar) \leq  [\tilde{G}_p(\ball)]^2,$$  with equality if and only if $K$ is an origin-symmetric ellipsoid.\\
\noindent (ii). For $-n\neq p<0$, $$\tilde{G}_p(K)\tilde{G}_p(\polar) \geq  c^n [\tilde{G}_p(\ball)]^2,$$ 
 with $c$ the universal constant from the Bourgain-Milman inverse Santal\'{o} inequality \cite{BM}.}

 \vskip 2mm  We will prove the following affine isoperimetric inequality for $\tilde{G}_p(K)$.  See \cite{Lu1, Petty1974, Petty1985} for $p\geq 1$. 

\vskip 2mm \noindent {\bf Theorem \ref{isoperimetric:geominimal}}  {\em 
Let $K$ be a convex body with centroid (or the Santal\'{o} point) at the origin. 

\vskip 2mm \noindent (i). If $p\geq 0$, then
\begin{eqnarray*}
\frac{\tilde{G}_p(K)}{\tilde{G}_p(B^n_2)}\leq \min\left\{\bigg(\frac{|K|}{|\ball|}\bigg)^{\frac{n-p}{n+p}}, \ \bigg(\frac{|\polar|}{|B^n_2|}\bigg)^{\frac{p-n}{n+p}}\right\}.
\end{eqnarray*}
Equality holds for $p>0$ if and only if $K$ is an origin-symmetric ellipsoid. For $p=0$, equality holds trivially for all $K$.
\vskip 2mm \noindent
(ii). If
$-n<p<0$, then
\begin{eqnarray*}
\frac{\tilde{G}_p(K)}{\tilde{G}_p(B^n_2)}\geq
\left(\frac{|K|}{|B^n_2|}\right)^{\frac{n-p}{n+p}},
\end{eqnarray*} with
equality if and only if $K$ is an origin-symmetric ellipsoid.
\vskip 2mm \noindent 
   (iii). If $p < -n$, then \begin{eqnarray*}
\frac{\tilde{G}_p(K)}{\tilde{G}_p(B^n_2)}\geq
\left(\frac{|\polar|}{|B^n_2|}\right)^{\frac{p-n}{n+p}},
\end{eqnarray*} with
equality if and only if $K$ is an origin-symmetric ellipsoid. Moreover,
\begin{equation*}
\frac{\tilde{G}_p(K)}{\tilde{G}_p(B^n_2 )}\geq
c^{\frac{np}{n+p}}\left(\frac{|K|}{|\ball|}\right)^{\frac{n-p}{n+p}},
\end{equation*}
where $c$ is the constant in the Bourgain-Milman inverse Santal\'o inequality \cite{BM}. }

 \vskip 2mm We will show the monotonicity of $\tilde{G}_p(\cdot)$. This result will be used to obtain a stronger version of Theorem \ref{isoperimetric:geominimal}. 
\vskip 2mm \noindent  {\bf Theorem \ref{monotone:geominimal:p--1}} {\em Let $K $ be a convex body with the origin in its interior and $p, q\neq 0$. \\ \noindent (i). If either $-n<q<p$ or $q<p<-n$, one has $$\bigg( \frac{\tilde{G}_q(K)}{n|K|}\bigg)^{\frac{n+q}{q}}\leq \bigg( \frac{\tilde{G}_p(K)}{n|K|}\bigg)^{\frac{n+p}{p}}.$$  \noindent (ii). If   $q<-n<p$, one has $$\bigg( \frac{\tilde{G}_q(K)}{n|K|}\bigg)^{\frac{n+q}{q}}\geq \bigg( \frac{\tilde{G}_p(K)}{n|K|}\bigg)^{\frac{n+p}{p}}.$$}
 
  This paper is organized as follows. Basic background and notation for convex geometry are given in Section \ref{section 2}.  In Section \ref{section 3}, we will provide our definition for the $L_p$ geominimal surface area of $K$. Such a definition is motivated by Theorem \ref{equivalent:affine:surface:area}. We will prove some properties of the $L_p$ geominimal surface area, such as affine invariance. Important inequalities, e.g., the affine isoperimetric inequality (i.e., Theorem \ref{isoperimetric:geominimal})  and a Santal\'{o} style inequality (i.e., Theorem \ref{geominimal product:bound}) are established in Section \ref{section 4}.  Comparison between the $L_p$ geominimal surface area and the $p$-surface area is also provided in Section \ref{section 4}. Cyclic inequalities and monotonicity properties for the  $L_p$ geominimal surface area are proved in Section \ref{section 5}. A stronger version of the affine isoperimetric inequality is also established in Section \ref{section 5}. 
  
  General references for convex geometry are \cite{Bus1958, Gruber2007, Lei1980, Sch}.

 
\section{Background and Notation}\label{section 2}
Our setting is $\bbR^n$ with the standard inner product   $\langle
\cdot, \cdot\rangle$, which induces the Euclidian norm $\|\cdot\|$. We use $\ball=\{x\in \bbR^n: \|x\|\leq 1\}$ to denote the unit Euclidean ball in $\bbR^n$ and $S^{n-1}=\{x\in \bbR^n: \|x\|=1\}$ to denote the unit sphere in $\bbR^n$ (i.e.,  the boundary of $B^n_2$). The volume of $\ball$ is denoted by $\o_n=|\ball|$. A {\it convex body} $K\subset \bbR^n$ is a convex compact subset of
$\bbR^n$ with nonempty interior. We use $\cK$ to denote the set of all convex bodies in $\bbR^n$. Through this paper, we will concentrate on the subset $\cK_0$ of $\cK$, which contains all convex bodies with the origin in their interior. The subset $\cK_c$ of $\cK_0$ contains all convex bodies in $\cK_0$ with centroid at the origin. The boundary of $K$ is denoted by $\partial K$. The
{\it polar body} of $K$, denoted by $\polar$, is defined as
$$\polar=\{y\in \bbR^n: \langle x,y \rangle \leq 1, \forall x\in
K\}.$$  For any convex body $K\in \cK_0$, its polar body $\polar\in \cK_0$ as well. The bipolar theorem (see
\cite{Sch}) states that $(K^\circ)^\circ =K$.  We write $T$ for an invertible linear transform from $\bbR^n$ to $\bbR^n$. The absolute value of the determinant of $T$ is written as $|det(T)|$.  A convex body $\E\in \cK$ is said to be an {\it origin-symmetric ellipsoid} if $\E=T\ball$ for some invertible linear transform $T$ on $\bbR^n$.

\vskip 2mm The {\it support
function} of $K\in \cK_0$ at the direction $u\in S^{n-1}$ is defined as
$$h_K(u)=\max _{x\in K} \langle x,u \rangle .$$
The support function $h_K(\cdot)$ determines a convex body uniquely, and can be used to calculate the volume of $\polar$, namely,  $$|\polar|= \frac{1}{n}\int _{S^{n-1}}
\frac{1} {h _K^n(u)}\,d\s(u), $$  where  $\s$
is the usual spherical measure on $S^{n-1}$. More generally,
for a set $M$, we use $|M|$ to denote the Hausdorff content of the
appropriate dimension. The Hausdorff metric $d_H$ is a natural metric for $\cK_0$. For two convex bodies $K, K'\in \cK_0$, the Hausdorff distance between $K$ and $K'$ is $$d_H(K, K')=\|h_K-h_{K'}\|_{\infty}=\sup_{u\in S^{n-1}}|h_K(u)-h_{K'}(u)|. $$

 For two convex bodies $K, L\in \cK_0$ and  $\l, \eta \geq 0$ (not both zero), the
{\it Minkowski linear combination} $\l K+ \eta L$ is the convex body
with support function $h_{\l K+\eta L}$, such that, 
$$h_{\l K+\eta L}(u)=\lambda h_K(u)+\eta h_L(u), \ \ \  \forall u\in S^{n-1}.$$  The {\it mixed volume} of $K$ and $L$, denoted by $V_1
(K,L)$, is defined by
$$V_1 (K,L)=\lim _{\e \rightarrow 0} \frac{|K+\e L|-|K|}{n\e}. $$ For each convex body $K\in \cK$, 
there is a positive Borel measure $S(K, \cdot)$ on $S^{n-1}$ (see \cite{Ale1937-1, Fenchel1938}), such
that, for all convex bodies $L$, $$V_1 (K,L)=\frac{1}{n}\int
_{S^{n-1}}h_L(u)\,dS(K,u).$$ As noted in \cite{Lu2}, the measure
$S(K,\cdot)$ has the following simple geometric description: for
any Borel subset $A$ of $S^{n-1}$, 
$$S(K,A)=|\{x\in \partial K:\ \mbox{$\exists u\in A$, s.t., $H(x,u)$ is a support hyperplane of $\partial K$ at $x$}
\}|.$$  If the measure $S(K, \cdot)$ is absolute continuous with respect to the spherical measure $\s$, then by the Radon-Nikodym theorem, there is a function $f_K: S^{n-1}\rightarrow \bbR$, the {\it curvature function} of $K$, such that, $$dS(K,u)=f_K(u)d\s(u).$$
 Let $\cF\subset \cK$ be the set of all convex bodies with curvature function. Respectively, we denote $\cF_0=\cF\cap \cK_0$  and $\cF_c=\cF\cap \cK_c$. We write $K\in \cF_0^+$ for those convex bodies $K\in \cF_0$  with continuous and positive curvature function $f_K(\cdot)$  on $S^{n-1}$.

\vskip 2mm \par The Minkowski linear combination and the mixed
volume $V_1(\cdot,\cdot)$ can be generalized to all $p\geq 1$. For
convex bodies $K, L\in \cK_0$, and $\l, \eta \geq 0$ (not both zeros), the {\it Firey
$p$-sum} $\l K+_p \eta L$ for $p\geq 1$ \cite{Firey1962} is the
convex body with support function defined by
$$\big(h_{\l K+_p\eta L}(u)\big)^p=\lambda
\big(h_K(u)\big)^p+\eta \big(h_L(u)\big)^p.$$ The {\it $p$-mixed volume}  of $K$ and $L$, denoted by
$V_p(K,L)$, was defined by Lutwak \cite{Lut1993} as \begin{equation*}
\frac{1}{p}\ V_p(K,L)=\lim _{\e
\rightarrow 0} \frac{|K+_p\e L|-|K|}{n\e}. \end{equation*} Lutwak \cite{Lut1993} proved that, for
$K\in \cK_0$, there is a measure $S_p(K,\cdot)$,
such that
\begin{equation} V_p(K,L)=\frac{1}{n}\int _{S^{n-1}}h_L^p(u)\,dS_p(K,u),\label{mixed:p:volume} \end{equation} holds for all $L\in \cK_0$.  The measure $S_p(K,\cdot)$ is absolutely continuous
with respect to the measure $S(K, \cdot)$. Moreover, for all $p\geq 1$,
\be\,dS_p(K,u)=h_K^{1-p}(u)\,dS(K, u). \label{p:surface:area}\ee Therefore, for all $K\in \cF_0$ and all $p\geq 1$, one can define the {\it $L_p$ curvature function} for $K$, denoted by $f_p(K, u)$, to be (see \cite{Lu1})  $$f_p(K, u)=h_K(u)^{1-p}f_K(u).$$

\vskip 2mm We will use $\cK_s\subset \cK_0$  to denote the subset of $\cK_0$ containing all convex bodies with Santal\'{o} point at the origin.  Here, a convex body $K$ is said to have Santal\'{o} point at the origin, if its polar body $\polar$ has centroid at the origin, i.e., $K\in \cK_s \Leftrightarrow \polar \in \cK_c$. We also denote by $\cF_s=\cF_0\cap \cK_s$ the subset of $\cF_0$ containing all convex bodies with curvature function and Santal\'{o} point at the origin.

A subset $L\subset \bbR^n$ is {\it star-shaped} (about $x_0 \in L$)
if there exists $x_0\in L$, such that, the line segment 
from $x_0$ to any point $x\in L$ is contained in $L$. Hereafter, we only work on $\cS_0$, the set of all $n$-dimensional {\it star bodies} about the origin (i.e., compact star-shaped subsets of $\bbR^n$ about the origin with
continuous and positive radial functions). The {\it radial function}
of $L\in \cS_0$ at the direction $u\in S^{n-1}$ is given by $$\r _L(u)=\max \{\l: \l u\in L\}.$$ Note that, the volume of $L$ can be calculated by $$|L|=\frac{1}{n}\int_{S^{n-1}}\r_L(u)^n\,d\s(u).$$
  For $K\in \cK_0$, $L\in \cS_0$, and $p\geq 1$,
let $V_p(K, L^\circ)$ be $$ V_p(K, L^\circ)=\frac{1}{n} \int
_{S^{n-1}} \r _L(u)^{-p} \,dS_p(K,u).$$ When $L\in \cK_0$, this is consistent with formula  (\ref{mixed:p:volume}) because $\r_L(u){h_{L^\circ}(u)}=1$,  $\forall u\in S^{n-1}$.

\section{$L_p$ geominimal surface area}\label{section 3}

\vskip 2mm \par For $p\geq 1$ and $K\in \cK_0$, Lutwak in \cite{Lu1} defined the $L_p$ geominimal surface area of $K$, denoted by $G_p(K)$, as \be  \o_n^{p/n}G_p(K)=\inf_{Q\in \cK_0}\{n V_p(K, Q)|Q^\circ|^{p/n} \}. \nonumber \ee 

We now propose a definition for the $L_p$ geominimal surface area of a convex body $K\in \cK_0$, denoted by $\tilde{G}_p(K)$, for all $-n\neq p\in \bbR$.  Before we state our definition, we need some notation. First, notice that the measure $S_p(K, \cdot)$ in (\ref{p:surface:area}), the $p$-mixed volume in (\ref{mixed:p:volume}), and the $L_p$ curvature function $f_p(K, \cdot)$  were only defined for $p\geq 1$ in \cite{Lu1}. However, we  can extend them to all $ p\in \bbR$. For two convex bodies $K, Q\in \cK_0$, we  define the $p$-mixed volume of $K$ and $Q$ by $$ V_p(K, Q)=\frac{1}{n} \int
_{S^{n-1}} h_Q(u)^p \,dS_p(K,u), \ \ \ p\in \bbR, $$ where $S_p(K,\cdot)$ is the $p$-surface area measure of $K\in \cK_0$  defined as $$\,dS_p(K,u)=h_K^{1-p}(u)\,dS(K, u),  \ \ \ p\in \bbR.$$  We define the {\it $L_p$ curvature function} (and denoted by $f_p(K, \cdot)$) for $K\in \cF_0$ by  $$f_p(K, u)=h_K(u)^{1-p}f_K(u), \ \ p\in \bbR. $$  Note that, if $K\in \cF_0$, then $$\,dS_p(K,u)=f_p(K, u)\,d\s(u),  \ \ \ p\in \bbR.$$ For $K\in \cK_0$ and  $L\in \cS_0$, let $V_p(K, L^\circ)$ be $$ V_p(K, L^\circ)=\frac{1}{n} \int
_{S^{n-1}} \r _L(u)^{-p} \,dS_p(K,u), \ \ \  p\in \bbR.$$
\bd Let $K\in \cK_0$ be a convex body with the origin in its interior. \\
  (i). For $p\geq 0$, we define the $L_p$ geominimal surface area of $K$ by \be\label{Lp:geominimal:Surface:Area:positive}  \tilde{G}_p(K)=\inf_{Q\in \cK_0} \left\{n  V_p(K, Q) ^{\frac{n}{n+p}}\
  |Q^\circ|^{\frac{p}{n+p}}\right\}.\ee
  (ii). For $-n\neq p<0$, we define the $L_p$ geominimal surface area of $K$ by \be\label{Lp:geominimal:Surface:Area:negative} \tilde{G}_p(K)=\sup_{Q\in \cK_0} \left\{n  V_p(K, Q) ^{\frac{n}{n+p}}\
    |Q^\circ|^{\frac{p}{n+p}}\right\}.\ee  \ed  
\noindent {\bf Remark.} Let $p=0$ and $Q\in \cK_0$ be any fixed convex body, one has $V_p(K, Q)=|K|$ and $ n  V_p(K, Q) ^{\frac{n}{n+p}} 
  |Q^\circ|^{\frac{p}{n+p}} =n|K|$. Therefore, $$\tilde{G}_p(K)=\inf_{Q\in \cK_0} \left\{n V_p(K, Q) ^{\frac{n}{n+p}}
  |Q^\circ|^{\frac{p}{n+p}}\right\}=n|K|.$$ The case $p=-n$ is not covered mainly because $n+p=0$ appears in the denominators of $\frac{p}{n+p}$ and $\frac{n}{n+p}$. More general $L_p$ geominimal surface area $\tilde{G}_p(K, \cM)$ can be defined with the infimum or the supremum taking over $\cM\subset \cK_0$.  This paper only deals with $\tilde{G}_p(K)= \tilde{G}_p(K, \cK_0)$. 
 For $p\geq 1$, our $L_p$ geominimal surface area is related to Lutwak's  by the following formula: for any $K\in \cK_0$,  \be [\tilde{G}_p(K)]^{n+p}=(n\o_n)^p [G_p(K)]^n.  \label{relation tilde G}\ee

 \bp\label{homogeneous:degree}  Let $K\in \cK_0$. For all $p\neq -n$ and all invertible linear transformation $T:\bbR^n\rightarrow \bbR^n$, one has $$ \tilde{G}_p(TK)=|det(T)|^{\frac{n-p}{n+p}} \tilde{G}_p(K).$$ \ep 
\noindent {\bf Proof.}  First, we note that for all $y\in TK$, there exists a unique $x\in K$, such that $y=Tx$ (as $T$ is invertible). The support function of $TK$ can be calculated as follows $$h_{TK}(v)=\max_{y\in TK} \langle y, v\rangle=\max_{x\in K} \langle Tx, v\rangle =\max_{x\in K} \langle x, T^{*} v\rangle=\|T^*v\| \max_{x\in K} \langle x, u \rangle=\|T^*v\| h_K(u),$$
where $T^{*}$ is the transpose of $T$ and $u=\frac{T^{*} v}{\|T^{*} v\|} $. Hence, $$\frac{h_{TQ}(v)}{h_{TK}(v)}=\frac{h_{Q}(u)}{h_{K}(u)}.$$ On the other hand, $\frac{1}{n} h_K(u)\,dS(K, u)$ is the volume element of $K$ and hence, $$h_{TK}(v)\,dS(TK, v)=|det(T)| h_K(u)\,dS(K, u).$$ Therefore, one has 
\begin{eqnarray*}
nV_p(TK, TQ)&=&\int _{S^{n-1}} \bigg(\frac{h_{TQ}(v)}{h_{TK}(v)}\bigg)^ph_{TK}(v)\,dS(TK, v)\\&=&|det(T)|  \int _{S^{n-1}} \bigg(\frac{h_{Q}(u)}{h_{K}(u)}\bigg)^ph_{K}(u)\,dS(K, u) \\&=&n|det(T)| V_p(K, Q). \end{eqnarray*} Recall that $(TQ)^\circ=[T^*]^{-1} Q^\circ$. This further implies that, for $p\geq 0$, 
\begin{eqnarray*}
\tilde{G}_p(TK)&=&\inf _{Q\in \cK_0} \{n [V_p(TK, TQ)]^{\frac{n}{n+p}} |(TQ)^\circ|^{\frac{p}{n+p}}\}\\ &=& |det(T)| ^{\frac{n}{n+p}}  |det([T^*]^{-1})|^{\frac{p}{n+p}} \inf _{Q\in \cK_0} \{n [V_p(K, Q)]^{\frac{n}{n+p}} |Q^\circ|^{\frac{p}{n+p}}\}\\&=&|det(T)|^{\frac{n-p}{n+p}} \tilde{G}_p(K).
\end{eqnarray*} Similarly, for $-n\neq p<0$, one has \begin{eqnarray*}
\tilde{G}_p(TK)&=&\sup _{Q\in \cK_0} \{n [V_p(TK, TQ)]^{\frac{n}{n+p}} |(TQ)^\circ|^{\frac{p}{n+p}}\}\\ &=& |det(T)| ^{\frac{n}{n+p}}  |det([T^*]^{-1})|^{\frac{p}{n+p}} \sup _{Q\in \cK_0} \{n [V_p(K, Q)]^{\frac{n}{n+p}} |Q^\circ|^{\frac{p}{n+p}}\}\\&=&|det(T)|^{\frac{n-p}{n+p}} \tilde{G}_p(K).
\end{eqnarray*}

\noindent {\bf Remark.} The $L_p$ geominimal surface area is invariant under all invertible linear transformations $T$ with $|det(T)|=1$, i.e., $ \tilde{G}_p(K)= \tilde{G}_p(TK)$.  Moreover, if $T$ is just a dilation, say $TK=rK$ with $r>0$, then, $$ \tilde{G}_p(rK)=r^{\frac{n (n-p)}{n+p}} \tilde{G}_p(K).$$ Note that the classical geominimal surface area $\tilde{G}_1(\cdot)$ is also translation invariant. That is, for all $K\in \cK_0$ and all $z_0\in \bbR^n$, one has $$\tilde{G}_1(K-z_0)=\tilde{G}_1(K). $$ This can be easily seen from the translation invariance of the surface area measure $\,dS(K, \cdot)$, which clearly implies, $\forall Q\in \cK_0$,  $$nV_1(K, Q)=\int_{S^{n-1}} h_Q \,dS(K, u)=\int_{S^{n-1}} h_Q \,dS(K-z_0, u)=nV_1(K-z_0, Q).$$
However, one {\it cannot} expect the translation invariance for $\tilde{G}_p(\cdot)$ for $-n\neq p\in \bbR$ (expect $p=1$ and $p=0$). As an example which will be used in Section \ref{section 5}, we show the following result. 
\bp\label{translation:invariance:fail:ball} Let $\E\subset \bbR^n$ be an origin-symmetric ellipsoid, and $z_0\neq 0$ be any interior point of $\E$.  
\\ (i). For all $p\in (0, 1)$, $$\tilde{G}_p(\E-z_0) <\tilde{G}_p(\E).$$ (ii). For all $p\in (-n, 0)$, $$\tilde{G}_p(\E-z_0) >\tilde{G}_p(\E).$$
\ep
\noindent {\bf Proof.}  Without loss of generality, one only needs to verify Proposition \ref{translation:invariance:fail:ball} for $\E=\ball$ due to  Proposition \ref{homogeneous:degree}.  In this case,  the nonzero interior point $z_0\in \ball$ satisfies $0<\|z_0\|<1$. Denote by $B_{z_0}=\ball-z_0$ the ball with center $z_0$ and radius $1$. 

\vskip 2mm \noindent
(i). For $p\in (0, 1)$, the function $g(t)=t^{1-p}$ is strictly concave on $t\in (0, \infty)$. By Jensen's inequality, one has, \begin{eqnarray*}
\frac{V_p(B_{z_0}, \ball)}{|\ball|} = \frac{1}{n|\ball|} \int_{S^{n-1}} [h_{B_{z_0}}(u)]^{1-p}\,d\s  \leq\!\left(\! \frac{1}{n|\ball|} \int_{S^{n-1}} h_{B_{z_0}}(u)\,d\s\! \right)^{1-p}\!\!=\!\!\left(\!\frac{|B_{z_0}|}{|\ball|}\!\right)^{1-p}=1.\end{eqnarray*} Equality holds if and only if $h_{B_{z_0}}(u)$ is a constant, which is not possible as $z_0\neq 0$. Therefore, as $p\in (0, 1)$, \begin{eqnarray*}\tilde{G}_p(B_{z_0})&\leq& n V_p(B_{z_0}, \ball)^{\frac{n}{n+p}} |\ball|^{\frac{p}{n+p}}<  n   |\ball| =\tilde{G}_p(\ball),\end{eqnarray*} where $\tilde{G}_p(\ball)=n|\ball|$ will be proved in formula (\ref{geominimal:ball}). 

\vskip 2mm\noindent (ii). For $p\in (-n, 0)$,  the function $g(t)=t^{1-p}$ is strictly convex on $(0, \infty)$. By Jensen's inequality, one has, \begin{eqnarray*}
\frac{V_p(B_{z_0}, \ball)}{|\ball|} =\frac{1}{n|\ball|} \int_{S^{n-1}} [h_{B_{z_0}}(u)]^{1-p}\,d\s  \geq \left(\! \frac{1}{n|\ball|} \int_{S^{n-1}} h_{B_{z_0}}(u)\,d\s \right)^{1-p}\!\!=\!\!\left(\!\frac{|B_{z_0}|}{|\ball|}\!\right)^{1-p}\!\!=1.\end{eqnarray*} Equality holds if and only if $h_{B_{z_0}}(u)$ is a constant, which is not possible as $z_0\neq 0$. Therefore, as $p\in (-n, 0)$, \begin{eqnarray*}\tilde{G}_p(B_{z_0})&\geq& n V_p(B_{z_0}, \ball)^{\frac{n}{n+p}} |\ball|^{\frac{p}{n+p}}> n   |\ball| =\tilde{G}_p(\ball),\end{eqnarray*} where $\tilde{G}_p(\ball)=n|\ball|$ is given by formula (\ref{geominimal:ball}). 

\vskip 2mm 
  Our definition of $\tilde{G}_p(K)$ is motivated by Theorem \ref{equivalent:affine:surface:area} regarding a recent extension of the $L_p$ affine surface area. Recall that, Lutwak in \cite{Lu1} defined the $L_p$ affine surface area of $K\in \cK_0$ for $p\geq 1$ by 
\begin{equation}\label{Lp-affine-surface-area-Lutwak}
as_p(K)=\inf_{L\in \cS_0} \left\{n\ V_p(K, L^\circ) ^{\frac{n}{n+p}}\
|L|^{\frac{p}{n+p}}\right\}.
\end{equation}  This definition generalizes the definition of the classical
affine surface area (for $p=1$) in \cite{Lei1989}. It was also proved  that, for all $K\in \cF_0^+$, 
$as_p(K)$ for $p\geq 1$ defined in
(\ref{Lp-affine-surface-area-Lutwak}) can be written as an
integral over $S^{n-1}$ (see Theorem 4.4 in \cite{Lu1}),
\be
as_{p}(K)=\int_{S^{n-1}} {f_p({K}, u)^{\frac{n}{n+p}}} d\sigma(u).  \label{Lp-affine-surface-area-integral}
\ee   Such an integral expression for $as_p(K)$ has been used to extend the $L_p$ affine surface area from $p\geq 1$ to all $-n\neq p\in \bbR$ (see e.g. \cite{MW2, SW5, SW4}).

\bt\label{equivalent:affine:surface:area}  Let $K\in \cF_0^+$ be a convex body with continuous and positive curvature function and with the origin in its interior.   \\ (i). For $p\geq 0$,  one has \begin{equation*} as_p(K)=\inf_{ L\in \cS_0} \left\{n V_p(K, L^\circ) ^{\frac{n}{n+p}}
|L|^{\frac{p}{n+p}}\right\}.
\end{equation*}  (ii). For $-n\neq p<0$, one has \begin{equation*} as_p(K)=\sup_{L\in \cS_0} \left\{nV_p(K, L^\circ) ^{\frac{n}{n+p}}
|L|^{\frac{p}{n+p}}\right\}.
\end{equation*} \et

 \vskip 2mm\noindent {\bf Proof.} (i). It is clear that if $p=0$, then $V_p(K, L^\circ)= |K|$, and hence $$as_0(K)=n|K|= \inf_{L\in \cS_0} \left\{n V_p(K, L^\circ) ^{\frac{n}{n+p}}
|L|^{\frac{p}{n+p}}\right\}. $$ Let $p\in (0, \infty)$. Recall that, for all $K\in \cF_0^+$, the $L_p$ affine surface area of $K$ is given by formula (\ref{Lp-affine-surface-area-integral}) 
$$as_{p}(K)=\int_{S^{n-1}} {f_p({K}, u)^{\frac{n}{n+p}}} d\sigma(u).$$ Therefore, for all   $L\in \cS_0$, one has,   \begin{eqnarray}as_{p}(K) &=&\int_{S^{n-1}} [\r_L ^{-p}(u){f_p({K}, u)]^{\frac{n}{n+p}}} \r_L(u)^{\frac{pn}{n+p}} d\sigma(u)\nonumber\\ 
&\leq & \left(\int_{S^{n-1}} \r_L ^{-p}(u)f_p({K}, u)d\sigma(u)\right) ^{\frac{n}{n+p}} \left(\int_{S^{n-1}}   \r_L(u)^{n} d\sigma(u)\right)^{\frac{p}{n+p}}\nonumber \\ &=& n [V_p(K, L^\circ)]^{\frac{n}{n+p}} |L|^{\frac{p}{n+p}}, \label{affine:Holder} 
\end{eqnarray} where the inequality follows from the H\"{o}lder inequality (see \cite{HLP}) and $p\in (0, \infty)$ (which implies $\frac{n}{n+p}\in (0, 1)$).   Taking the infimum over  $L\in \cS_0$, one gets, for all $p\in (0, \infty)$,  
 \begin{equation*} as_p(K)\leq \inf_{ L\in \cS_0} \left\{n V_p(K, L^\circ) ^{\frac{n}{n+p}}
|L|^{\frac{p}{n+p}}\right\}. 
\end{equation*}  
Note that $K\in \cF_0^+$ and hence $f_K(u)>0$ for all $u\in S^{n-1}$. Equality holds in (\ref{affine:Holder}) if and only if, for  some constant $a>0$, $$a^{n+p} \r_L ^{-p}(u){f_p({K}, u)}= \r_L(u)^{n} \Leftrightarrow a^{n+p} f_p(K, u)=\r_L(u)^{n+p}, \ \ \forall u\in S^{n-1}. $$  Thus, one can define the star body $L_0\in \cS_0$ by the radial function, which is clearly continuous and positive,  $$\r _{L_0}(u)=a [f_p(K, u)]^{\frac{1}{n+p}}>0, \ \ \forall u\in S^{n-1},$$ and hence $$as_p(K)=n [V_p(K, L_0^\circ)]^{\frac{n}{n+p}} |L_0|^{\frac{p}{n+p}}=\inf_{L\in \cS_0} \{n [V_p(K, L^\circ)]^{\frac{n}{n+p}} |L|^{\frac{p}{n+p}}\}. $$

\vskip 2mm \noindent (ii).   Let $-n\neq p\in (-\infty, 0)$. Similar to (\ref{affine:Holder}), for all   $L\in \cS_0$, one has,   \begin{eqnarray}as_{p}(K) &=&\int_{S^{n-1}} [\r_L ^{-p}(u){f_p({K}, u)]^{\frac{n}{n+p}}} \r_L(u)^{\frac{pn}{n+p}} d\sigma(u)\nonumber\\ 
&\geq & \left(\int_{S^{n-1}} \r_L ^{-p}(u)f_p({K}, u)d\sigma(u)\right) ^{\frac{n}{n+p}} \left(\int_{S^{n-1}}   \r_L(u)^{n} d\sigma(u)\right)^{\frac{p}{n+p}}\nonumber \\ &=& n [V_p(K, L^\circ)]^{\frac{n}{n+p}} |L|^{\frac{p}{n+p}}, \label{affine:Holder:negative} 
\end{eqnarray} where the inequality follows from the H\"{o}lder inequality (see \cite{HLP}) and $-n\neq p\in (-\infty, 0)$ (which implies either $\frac{n}{n+p}>1$ or $\frac{n}{n+p}<0$).   Taking the supremum over 
  $L\in \cS_0$, one gets, for all $-n\neq p\in (-\infty, 0)$,  
 \begin{equation*} as_p(K)\geq \sup_{L\in \cS_0} \left\{n V_p(K, L^\circ) ^{\frac{n}{n+p}}
|L|^{\frac{p}{n+p}}\right\}. 
\end{equation*} 
Note that $K\in \cF_0^+$ and hence $f_K(u)>0$ for all $u\in S^{n-1}$. Equality holds in (\ref{affine:Holder:negative}) if and only if, for some constant $a>0$,  $$a^{n+p} \r_L ^{-p}(u){f_p({K}, u)}= \r_L(u)^{n} \Leftrightarrow a^{n+p} f_p(K, u)=\r_L(u)^{n+p}, \ \ \forall u\in S^{n-1}. $$  Thus, one can define the star body $L_0\in \cS_0$ by the radial function, which is clearly continuous and positive,  $$\r _{L_0}(u)=a [f_p(K, u)]^{\frac{1}{n+p}}>0, \ \ \forall u\in S^{n-1},$$ and hence $$as_p(K)=n [V_p(K, L_0^\circ)]^{\frac{n}{n+p}} |L_0|^{\frac{p}{n+p}}=\sup_{L\in \cS_0} \{n [V_p(K, L^\circ)]^{\frac{n}{n+p}} |L|^{\frac{p}{n+p}}\}. $$

\noindent {\bf Remark.} The proof of Theorem \ref{equivalent:affine:surface:area} implies that, for all $K\in \cF_0^+$ and all $-n\neq p\in \bbR$, \be \label{affine=geominimal} [as_p(K)]^{n+p}=n^{n+p} \o_n^n |\L _p K|^p,$$ where  $\L_p K\in \cS_0$ is the $p$-curvature image of $K$ defined by (for $p\geq 1$, see \cite{Lu1}) $$f_p(K, u)=\frac{\o_n}{|\L_pK|}[\r_{\L_pK}(u)]^{n+p}.\ee

Motivated by Theorem \ref{equivalent:affine:surface:area}, one may define the $L_p$ affine surface area for all $-n\neq p\in \bbR$ as follows.  \bd Let $K\in \cK_0$ be a convex body with the origin in its interior. \\
  (i). For $p=0$, let $as_p(K)=n|K|$. For $p>0$, let \be\label{Lp:surface:Surface:Area:positive---1}  as_p(K)=\inf_{L\in \cS_0} \left\{n V_p(K, L^\circ) ^{\frac{n}{n+p}}
  |L|^{\frac{p}{n+p}}\right\}.\ee
  (ii). For $-n\neq p<0$, let \be\label{Lp:affine:Surface:Area:negative---1} as_p(K)=\sup_{L\in \cS_0} \left\{n V_p(K, L^\circ) ^{\frac{n}{n+p}}
    |L|^{\frac{p}{n+p}}\right\}.\ee  \ed  
One can easily check: for all $K\in \cK_0$,  $as_p(K)\leq \tilde{G}_p(K)$ for $p\in (0, \infty);$ while $as_p(K)\geq \tilde{G}_p(K)$ for $-n\neq p\in (-\infty, 0)$. To see this, by formula (\ref{Lp:surface:Surface:Area:positive---1}), one has for $p>0$, \begin{eqnarray} as_p(K)&=&\inf_{ L\in \cS_0} \left\{n V_p(K, L^\circ) ^{\frac{n}{n+p}}
|L|^{\frac{p}{n+p}}\right\}\nonumber \\ &\leq&  \inf_{ L\in \cK_0} \left\{n V_p(K, L^\circ) ^{\frac{n}{n+p}}
|L|^{\frac{p}{n+p}}\right\} = \tilde{G}_p(K),\label{affine<geominimal} \end{eqnarray} where the inequality is due to $\cK_0\subset\cS_0$. By formula (\ref{Lp:affine:Surface:Area:negative---1}), one has for $-n\neq p<0$, \begin{eqnarray} as_p(K)&=&\sup_{ L\in \cS_0} \left\{n V_p(K, L^\circ) ^{\frac{n}{n+p}}
|L|^{\frac{p}{n+p}}\right\}\nonumber \\ &\geq & \sup_{ L\in \cK_0} \left\{n V_p(K, L^\circ) ^{\frac{n}{n+p}}
|L|^{\frac{p}{n+p}}\right\}= \tilde{G}_p(K),\label{affine>geominimal} \end{eqnarray} where the inequality is due to $\cK_0\subset \cS_0$.

\vskip 2mm 
Let $-n\neq p\in \bbR$. Define the subset $\cV_p$ of $\cF_0^+$ as 
$$\cV_p=\{K\in \cF^+_0: \ \exists Q\in \cK_0\ \ with\ \ f_p(K,u)=h_Q(u)^{-(n+p)}, \ \ \forall u\in S^{n-1}\}. $$ Clearly, $\cV_p\neq \emptyset$ as $\ball \in \cV_p$.   The following proposition describes the relation between $\L_p K$ and $\cV_p$. See \cite{Lu1} for $p\geq 1$. 
\bp\label{equivalent:p:image:curvature} Let $-n\neq p\in \bbR$ and $K\in \cF_0^+$, then $$K\in \cV_p \ \ \ \mbox{if and only if}\ \ \ \L_pK\in \cK_0. $$\ep
\noindent {\bf Proof.} Let $-n\neq p \in \bbR$, and $K\in \cF_0^+$, then \begin{eqnarray*}
K\in \cV_p &\Leftrightarrow& [f_p(K, u)]^{\frac{-1}{n+p}} =h_{Q}(u), \ \ \exists Q\in \cK_0,  \forall u\in S^{n-1}\\
&\Leftrightarrow& [f_p(K, u)]^{\frac{1}{n+p}} =\r_{Q^\circ}(u), \ \ \exists Q\in \cK_0,  \forall u\in S^{n-1}\\ &\Leftrightarrow& \left(\frac{\o_n}{|\L _p K|}\right)^{\frac{1}{n+p}}\r_{\L _p K}(u) =\r_{Q^\circ}(u), \ \ \exists Q\in \cK_0,  \forall u\in S^{n-1}\\  
 &\Leftrightarrow& \left(\frac{\o_n}{|\L _p K|}\right)^{\frac{1}{n+p}}{\L _p K} ={Q^\circ}\in \cK_0.
\end{eqnarray*}

\bp\label{geominimal:ball-1} Let $-n\neq p\in \bbR$ and $K\in \cV_p$, then  $\tilde{G}_p(K)=as_p(K).$ \ep
\noindent{\bf Proof.} The proposition holds for $p=0$ as $\tilde{G}_0(K)=as_0(K)=n|K|$ for all $K\in \cK_0$. Assume $p>0$. From Proposition \ref{equivalent:p:image:curvature}, $K\in \cV_p$ implies $\L_pK\in \cK_0$. Therefore, formula (\ref{affine=geominimal}) and inequality (\ref{affine<geominimal}) imply    \begin{eqnarray*}
 \tilde{G}_p(K)&\geq& as_p(K)=n \o_n^{\frac{n}{n+p}}|\L _pK|^{\frac{p}{n+p}}\\ &=& n V_p(K, (\L_pK)^\circ)^{\frac{n}{n+p}}|\L_pK|^{\frac{p}{n+p}}\\ &\geq& \inf_{Q\in \cK_0} \{n V_p(K, Q^\circ)^{\frac{n}{n+p}}|Q|^{\frac{p}{n+p}}\}\\ &=& \tilde{G}_p(K). 
\end{eqnarray*} Hence, $ \tilde{G}_p(K)= as_p(K)$ for all $K\in \cV_p$. 

\vskip 2mm\noindent  For $-n\neq p<0$,  formula (\ref{affine=geominimal}) and inequality (\ref{affine>geominimal}) imply   \begin{eqnarray*}
 \tilde{G}_p(K)&\leq& as_p(K)= n V_p(K, (\L_pK)^\circ)^{\frac{n}{n+p}}|\L_pK|^{\frac{p}{n+p}}\\ &\leq& \sup_{Q\in \cK_0} \{n V_p(K, Q^\circ)^{\frac{n}{n+p}}|Q|^{\frac{p}{n+p}}\}\\ &=& \tilde{G}_p(K). 
\end{eqnarray*} Hence, $ \tilde{G}_p(K)= as_p(K)$ for all $K\in \cV_p$.

\vskip 2mm \noindent {\bf Remark.} Recall that $\ball\in\cV_p$. Thus, \be \label{geominimal:ball}\tilde{G}_p(\ball)=as_p(\ball)=n\o_n.\ee  Hence, for all origin-symmetric ellipsoids $\E=T\ball$ with $T$ an invertible linear transform on $\bbR^n$, one has, by Proposition \ref{homogeneous:degree},  $$\tilde{G}_p(\E)=\tilde{G}_p(T\ball)=|det(T)|^{\frac{n-p}{n+p}} \tilde{G}_p(\ball)=|det(T)|^{\frac{n-p}{n+p}} n\o_n.$$

\section{Affine isoperimetric and Santal\'{o} style inequalities} \label{section 4} 
This section is mainly dedicated to the affine isoperimetric inequality and a Santal\'{o} style inequality for the $L_p$ geominimal surface area. 

\bp \label{bounded by volume product} Let $K\in \cK_0$ be a convex body with the origin in its interior. 
\\ (i). For $p\geq 0$, one has $$\tilde{G}_p(K) \leq n|K|^{\frac{n}{n+p}}|K^\circ|^{\frac{p}{n+p}}.$$ 
(ii). For $-n\neq p<0$, one has $$\tilde{G}_p(K) \geq n|K|^{\frac{n}{n+p}}|K^\circ|^{\frac{p}{n+p}}. $$  \ep
\noindent {\bf Proof.} 
 Note that $V_p(K, K)=|K|$ for all $K\in \cK_0$ and for all $-n\neq p\in \bbR$. 
 
 \vskip 2mm \noindent (i). The case $p=0$ is clear (and in fact ``=" always holds).  For $p>0$, by formula (\ref{Lp:geominimal:Surface:Area:positive}) and $K\in \cK_0$, one has 
 $$\tilde{G}_p(K)=\inf_{Q\in \cK_0} \left\{n V_p(K, Q) ^{\frac{n}{n+p}}
  |Q^\circ|^{\frac{p}{n+p}}\right\}\leq n|K|^{\frac{n}{n+p}}|K^\circ|^{\frac{p}{n+p}}.$$ 
 
\noindent (ii). For $-n\neq p<0$, by formula (\ref{Lp:geominimal:Surface:Area:negative}) and $K\in \cK_0$, one gets, $$\tilde{G}_p(K)=\sup_{Q\in \cK_0} \left\{n V_p(K, Q) ^{\frac{n}{n+p}}
  |Q^\circ|^{\frac{p}{n+p}}\right\}\geq n|K|^{\frac{n}{n+p}}|K^\circ|^{\frac{p}{n+p}}.$$

\bp \label{geominimal product}
 Let $K\in \cK_0$ be a convex body with the origin in its interior. 
\\ (i). For $p\geq 0$, one has $$\tilde{G}_p(K)\tilde{G}_p(\polar) \leq n^2 |K| |K^\circ| .$$  
(ii). For $-n\neq p<0$, one has $$\tilde{G}_p(K)\tilde{G}_p(\polar) \geq n^2 |K| |K^\circ|. $$  
\ep \noindent {\bf Proof.}  (i). For $p\geq 0$, using Proposition \ref{bounded by volume product} for both $K$ and $\polar$, one gets $$\tilde{G}_p(K)\tilde{G}_p(\polar) \leq (n|K|^{\frac{n}{n+p}}|K^\circ|^{\frac{p}{n+p}} )(n |K|^{\frac{p}{n+p}}|K^\circ|^{\frac{n}{n+p}})=n^2|K||\polar|.$$

\noindent (ii). For $-n\neq p<0$, using Proposition \ref{bounded by volume product} for both $K$ and $\polar$, one gets $$\tilde{G}_p(K)\tilde{G}_p(\polar) \geq (n|K|^{\frac{n}{n+p}}|K^\circ|^{\frac{p}{n+p}} )(n |K|^{\frac{p}{n+p}}|K^\circ|^{\frac{n}{n+p}})=n^2|K||\polar|.$$

The Blaschke-Santal\'{o} inequality provides a precise upper bound for $M(K)=|K||\polar|$ with $K\in \cK_c$ (or $K\in \cK_s$). That is, for all $K\in \cK_c$ (or $K\in \cK_s$), $$M(K)\leq M(\ball)=\o_n^2,$$ with equality if and only if $K$ is an origin-symmetric ellipsoid. Finding the precise lower bound for $M(K)$ is still an open problem and is known as the Mahler conjecture:  {the precise lower bound for $M(K)$ is conjectured to be obtained by the cube among all origin-symmetric convex bodies (i.e., $K=-K$),  and by the simplex among all convex bodies $K\in \cK_c$ (or $K\in \cK_s$).}  A remarkable result by Bourgain and Milman \cite{BM} states that, there is a universal constant $c>0$ (independent of $n$ and $K$), such that, for all $K\in \cK_c$ (or $K\in \cK_s$), \be M(K)\geq c^n M(\ball)=c^n\o_n^2. \nonumber \ee Nice estimates for the constant $c$ can be found in \cite{GK2, Nazarov2012}. 

\vskip 2mm We now prove the following Santal\'{o} style inequality for the $L_p$ geominimal surface area $\tilde{G}_p$. See \cite{Zhou2011} for $p\geq 1$. 
\bt \label{geominimal product:bound} Let $K\in \cK_c$ or $K\in \cK_s$. \\ (i). For $p\geq 0$,  $$\tilde{G}_p(K)\tilde{G}_p(\polar) \leq  [\tilde{G}_p(\ball)]^2,$$  with equality if and only if $K$ is an origin-symmetric ellipsoid.\\
\noindent (ii). For $-n\neq p<0$, $$\tilde{G}_p(K)\tilde{G}_p(\polar) \geq  c^n [\tilde{G}_p(\ball)]^2,$$ 
 with $c$ the universal constant from the Bourgain-Milman inverse Santal\'{o} inequality \cite{BM}.\et
 
 \noindent {\bf Proof.} (i). Recall that $\tilde{G}_p(\ball)=n|\ball|$. Proposition \ref{geominimal product} implies that \be \label{geominimal product:proof}\tilde{G}_p(K)\tilde{G}_p(\polar) \leq n^2|K||\polar|\leq n^2\o_n^2 =  [\tilde{G}_p(\ball)]^2,\ee where the second inequality is from the Blaschke-Santal\'{o} inequality. Clearly, Proposition \ref{homogeneous:degree} implies that equality holds in inequality (\ref{geominimal product:proof}) if $K$ is an origin-symmetric ellipsoid. On the other hand, the equality holds in inequality (\ref{geominimal product:proof}) only if the equality holds in the Blaschke-Santal\'{o} inequality, that is,  $K$ has to be an origin-symmetric ellipsoid. 
  
\vskip 2mm \noindent (ii).  Proposition \ref{geominimal product} implies $$\tilde{G}_p(K)\tilde{G}_p(\polar) \geq n^2|K||\polar|\geq c^n n^2\o_n^2=c^n  [\tilde{G}_p(\ball)]^2,$$ where the second inequality is from the Bourgain-Milman inverse  Santal\'{o} inequality. 

\vskip 2mm \noindent {\bf Remark.} Note that the Bourgain-Milman inverse  Santal\'{o} inequality still holds true for all $K\in \cK_0$. This is because $M(K)\geq M(K-z_0)\geq c^n \o_n^2$, where $z_0$ is the centroid of $K$ (and then $K-z_0\in \cK_c$). Hence,  part (ii) still holds for all $K\in \cK_0$. 

\vskip 2mm The related affine isoperimetric inequality for the $L_p$ geominimal surface area $\tilde{G}_p(K)$ is stated in the following theorem. 
\bt \label{isoperimetric:geominimal}
Let $K\in \cK_c$ or $K\in \cK_s$. 

\vskip 2mm \noindent (i). If $p\geq 0$, then
\begin{eqnarray*}
\frac{\tilde{G}_p(K)}{\tilde{G}_p(B^n_2)}\leq \min\left\{\bigg(\frac{|K|}{|\ball|}\bigg)^{\frac{n-p}{n+p}}, \ \bigg(\frac{|\polar|}{|B^n_2|}\bigg)^{\frac{p-n}{n+p}}\right\}.
\end{eqnarray*}
Equality holds for $p>0$ if and only if $K$ is an origin-symmetric ellipsoid. For $p=0$, equality holds trivially for all $K$.
\vskip 2mm \noindent
(ii). If
$-n<p<0$, then
\begin{eqnarray*}
\frac{\tilde{G}_p(K)}{\tilde{G}_p(B^n_2)}\geq
\left(\frac{|K|}{|B^n_2|}\right)^{\frac{n-p}{n+p}},
\end{eqnarray*} with
equality if and only if $K$ is an origin-symmetric ellipsoid.
\vskip 2mm \noindent 
(iii). If $p < -n$, then \begin{eqnarray*}
\frac{\tilde{G}_p(K)}{\tilde{G}_p(B^n_2)}\geq
\left(\frac{|\polar|}{|B^n_2|}\right)^{\frac{p-n}{n+p}},
\end{eqnarray*} with
equality if and only if $K$ is an origin-symmetric ellipsoid. Moreover, \begin{equation*}
\frac{\tilde{G}_p(K)}{\tilde{G}_p(B^n_2 )}\geq
c^{\frac{np}{n+p}}\left(\frac{|K|}{|\ball|}\right)^{\frac{n-p}{n+p}},
\end{equation*}
where $c$ is the constant in the Bourgain-Milman inverse Santal\'o inequality \cite{BM}.
\et
 \noindent {\bf Proof.} (i). The case $p=0$ is  trivial, and we only prove the case $p>0$. Combining
Proposition \ref{bounded by volume product}, the Blaschke-Santal\'o inequality, and
$\tilde{G}_p(B^n_2)=n|\ball|$,
one obtains
\begin{eqnarray*}
\frac{\tilde{G}_p(K)}{\tilde{G}_p(B^n_2)}\leq
\left(\frac{|K|}{|B^n_2|}\right)^{\frac{n}{n+p}}\left(\frac{|K^\circ|}{|B^n_2|}\right)^{\frac{p}{n+p}}= \left(\frac{|K|}{|B^n_2|}\right)^{\frac{n-p}{n+p}} \left(\frac{M(K)}{M(\ball)}\right)^{\frac{p}{n+p}}\leq
\left(\frac{|K|}{|B^n_2|}\right)^{\frac{n-p}{n+p}}.
\end{eqnarray*} Similarly, \begin{eqnarray*}
\frac{\tilde{G}_p(K)}{\tilde{G}_p(B^n_2)}\leq
\left(\frac{|K|}{|B^n_2|}\right)^{\frac{n}{n+p}}\left(\frac{|K^\circ|}{|B^n_2|}\right)^{\frac{p}{n+p}}= \left(\frac{|\polar|}{|B^n_2|}\right)^{\frac{p-n}{n+p}} \left(\frac{M(K)}{M(\ball)}\right)^{\frac{n}{n+p}}\leq
\left(\frac{|\polar|}{|B^n_2|}\right)^{\frac{p-n}{n+p}}.
\end{eqnarray*}
 Clearly, equality holds if $K$ is an origin-symmetric ellipsoid. On the other hand, equality holds in the above inequalities only if equality holds in the Blaschke-Santal\'{o} inequality, that is, $K$ has to be an origin-symmetric ellipsoid. 

\vskip 2mm\noindent  (ii). Let $-n<p<0$. Combining
Proposition \ref{bounded by volume product} and
$\tilde{G}_p(B^n_2)=n|\ball|$,
one obtains
\begin{eqnarray*}
\frac{\tilde{G}_p(K)}{\tilde{G}_p(B^n_2)}\geq
\left(\frac{|K|}{|B^n_2|}\right)^{\frac{n}{n+p}}\left(\frac{|K^\circ|}{|B^n_2|}\right)^{\frac{p}{n+p}}= \left(\frac{|K|}{|B^n_2|}\right)^{\frac{n-p}{n+p}} \left(\frac{M(K)}{M(\ball)}\right)^{\frac{p}{n+p}} \geq
\left(\frac{|K|}{|B^n_2|}\right)^{\frac{n-p}{n+p}}, 
\end{eqnarray*} where the last inequality follows from  the Blaschke-Santal\'o inequality and $\frac{p}{n+p}<0$. Clearly, equality holds if $K$ is an origin-symmetric ellipsoid. On the other hand, equality holds in the above inequalities only if equality holds in the Blaschke-Santal\'{o} inequality, that is, $K$ has to be an origin-symmetric ellipsoid.

\vskip 2mm \noindent (iii). Let $p<-n$. Combining
Proposition \ref{bounded by volume product} and
$\tilde{G}_p(B^n_2) =n|\ball|$,
one obtains
\begin{eqnarray*}
\frac{\tilde{G}_p(K)}{\tilde{G}_p(B^n_2)}\geq
\left(\frac{|K|}{|B^n_2|}\right)^{\frac{n}{n+p}}\left(\frac{|K^\circ|}{|B^n_2|}\right)^{\frac{p}{n+p}} =\left(\frac{|\polar|}{|B^n_2|}\right)^{\frac{p-n}{n+p}} \left(\frac{M(K)}{M(\ball)}\right)^{\frac{n}{n+p}} \geq  \left(\frac{|\polar|}{|B^n_2|}\right)^{\frac{p-n}{n+p}}, 
\end{eqnarray*} where the last inequality follows from  the Blaschke-Santal\'o inequality and $\frac{n}{n+p}<0$. Clearly, equality holds if $K$ is an origin-symmetric ellipsoid. On the other hand, equality holds in the above inequalities only if equality holds in the Blaschke-Santal\'{o} inequality, that is, $K$ has to be an origin-symmetric ellipsoid.

Moreover, from  $\frac{p}{n+p}>0$ and the Bourgain-Milman inverse Santal\'{o} inequality, one has
\begin{eqnarray*}
\frac{\tilde{G}_p(K)}{\tilde{G}_p(B^n_2)}\geq
\left(\frac{|K|}{|B^n_2|}\right)^{\frac{n}{n+p}}\left(\frac{|K^\circ|}{|B^n_2|}\right)^{\frac{p}{n+p}} = \left(\frac{|K|}{|B^n_2|}\right)^{\frac{n-p}{n+p}} \left(\frac{M(K)}{M(\ball)}\right)^{\frac{p}{n+p}} \geq c^{\frac{np}{n+p}}
\left(\frac{|K|}{|B^n_2|}\right)^{\frac{n-p}{n+p}}. 
\end{eqnarray*}

\vskip 2mm \noindent {\bf Remark.} If we assume that $|K|=|\ball|$, then part (i) of Theorem \ref{isoperimetric:geominimal} implies that $\tilde{G}_p(K)\leq \tilde{G}_p(\ball)$. Let $B_K$ be the origin-symmetric Euclidean ball with the same volume as $K$, and then the radius of $B_K$ is $r=\left(\frac{|K|}{|\ball|}\right)^{\frac{1}{n}}$.  Proposition \ref{homogeneous:degree} implies that $$\tilde{G}_p(B_K)=\left(\frac{|K|}{|B^n_2|}\right)^{\frac{n-p}{n+p}}\tilde{G}_p(\ball).$$ Therefore, part (i) of Theorem \ref{isoperimetric:geominimal} implies $$\tilde{G}_p(K)\leq \tilde{G}_p(B_K),  \ \ \ \ p>0,$$ with equality if and only if $K$ is an origin-symmetric ellipsoid. That is, {\em among all convex bodies in $\cK_c$ (or $\cK_s$) with fixed volume, the $L_p$ geominimal surface area for $p>0$ attains the maximum at and only at origin-symmetric ellipsoids.} Similarly, part (ii) of Theorem \ref{isoperimetric:geominimal} can be rewritten as $$\tilde{G}_p(K)\geq \tilde{G}_p(B_K), \ \ \ p\in(-n, 0),$$ with equality if and only if $K$ is an origin-symmetric ellipsoid. That is, {\em among all convex bodies in $\cK_c$ (or $\cK_s$) with fixed volume,   the $L_p$ geominimal surface area for $-n<p<0$ attains the minimum at and only at origin-symmetric ellipsoids.} From part (iii) of Theorem \ref{isoperimetric:geominimal}, one gets, for $K\in \cK_c$ (or $K\in \cK_s$),  $$\tilde{G}_p(K)\tilde{G}_p(B_{\polar})\geq [\tilde{G}_p(\ball)]^2, \ \ \ \ \ p<-n,$$  with equality if and only if $K$ is an origin-symmetric ellipsoid. From part (iii) of Theorem \ref{isoperimetric:geominimal} and the remark after Theorem \ref{geominimal product:bound},  one sees  also that  $$\tilde{G}_p(K)\geq c^{\frac{np}{n+p}}\tilde{G}_p(B_K), \ \ \ p<-n,$$ for all $K\in \cK_0$. For $p=1$, this is the classical result due to Petty \cite{Petty1974, Petty1985}, where the condition $K\in \cK_c$ or $K\in \cK_s$ can be replaced by $K\in \cK_0$ due to the translation invariant property of the classical geominimal surface area. The case $p> 1$ is due to Lutwak \cite{Lu1}.   A stronger version of Theorem \ref{isoperimetric:geominimal} for $p\in (-n, 1)$, where the condition on centroid or on Santal\'{o} point is removed, will be proved in Section \ref{section 5}.

\vskip 2mm From Lemma 2.1 in  \cite{Petty1974} (see page 79) and formula (\ref{relation tilde G}), one has $\tilde{G}_1(K_1)\leq \tilde{G}_1(K_2)$ for any two convex bodies $K_1, K_2\in \cK_0$ such that $K_1\subset K_2$. For the $L_p$ geominimal surface area, we have the following similar result.

\bc Let $\E$ be an origin-symmetric ellipsoid and $K\in \cK_c$ (or $K\in \cK_s$). \\ (i). For $p\in (0, n)$ and $K\subset \E$, one has $$\tilde{G}_p(K)\leq \tilde{G}_p(\E),$$ with equality if and only if $K=\E$. \\ (ii). For $p\in (n, \infty)$ and $ \E\subset K$, one has $$\tilde{G}_p(K)\leq \tilde{G}_p(\E),$$  with equality if and only if $K=\E$. \\  (iii). For $p\in (-n, 0)$ and $ \E\subset K$, one has $$\tilde{G}_p(K)\geq \tilde{G}_p(\E),$$ with equality if and only if $K=\E$. \\ (iv). For $p<-n$ and $ K\subset \E$,  one has $$\tilde{G}_p(K)\geq \tilde{G}_p(\E),$$ with equality if and only if $K=\E$.  \ec
\noindent {\bf Proof.} (i). Let $p\in (0, n)$ and hence $\frac{n-p}{n+p}>0$. Then it follows from $K\subset \E$ that $$\left(\frac{|K|}{|\E|}\right)^{\frac{n-p}{n+p}} \leq 1.$$ Combined with Theorem \ref{isoperimetric:geominimal} and Proposition \ref{homogeneous:degree}, one has  
\begin{eqnarray*}
\frac{\tilde{G}_p(K)}{\tilde{G}_p(\E)}\leq \left(\frac{|K|}{|\E|}\right)^{\frac{n-p}{n+p}} \leq 1\ \ \  \Rightarrow\ \ \  {\tilde{G}_p(K)}\leq {\tilde{G}_p(\E)}.
\end{eqnarray*}  Clearly equality holds if $K=\E$. On the other hand, if $K\subsetneq \E$, then $|K|<|\E|$ and hence ${\tilde{G}_p(K)}< {\tilde{G}_p(\E)}$. 

\vskip 2mm \noindent 
(ii). Let $p\in (n, \infty)$ and hence $\frac{n-p}{n+p}<0$. Then it follows from $\E\subset K$ that $$\left(\frac{|K|}{|\E|}\right)^{\frac{n-p}{n+p}} \leq 1.$$ Combined with Theorem \ref{isoperimetric:geominimal} and Proposition \ref{homogeneous:degree}, one has  
\begin{eqnarray*}
\frac{\tilde{G}_p(K)}{\tilde{G}_p(\E)}\leq \left(\frac{|K|}{|\E|}\right)^{\frac{n-p}{n+p}} \leq 1 \ \ \ \Rightarrow \ \ \ {\tilde{G}_p(K)}\leq {\tilde{G}_p(\E)}.
\end{eqnarray*}  Clearly equality holds if $K=\E$. On the other hand, if $\E\subsetneq K$, then $|\E|<|K|$ and hence ${\tilde{G}_p(K)}< {\tilde{G}_p(\E)}$. 

\vskip 2mm \noindent 
(iii). Let $p\in (-n, 0)$ and hence $\frac{n-p}{n+p}>0$. Then it follows from $\E\subset K$ that $$\left(\frac{|K|}{|\E|}\right)^{\frac{n-p}{n+p}} \geq 1.$$ Combined with Theorem \ref{isoperimetric:geominimal} and Proposition \ref{homogeneous:degree}, one has  
\begin{eqnarray*}
\frac{\tilde{G}_p(K)}{\tilde{G}_p(\E)}\geq \left(\frac{|K|}{|\E|}\right)^{\frac{n-p}{n+p}} \geq 1\ \ \  \Rightarrow\ \ \  {\tilde{G}_p(K)}\geq {\tilde{G}_p(\E)}.
\end{eqnarray*}  Clearly equality holds if $K=\E$. On the other hand, if $\E\subsetneq K$, then $|\E|<|K|$ and hence ${\tilde{G}_p(K)}> {\tilde{G}_p(\E)}$. 

\vskip 2mm \noindent 
(iv). Let $p<-n$ and hence $\frac{p-n}{n+p}>0$. Then it follows from $K\subset\E$ that $$\E^\circ\subset \polar \ \ \ \Rightarrow \ \ \ \left(\frac{|\polar|}{|\E^\circ|}\right)^{\frac{p-n}{n+p}} \geq 1.$$ Combined with Theorem \ref{isoperimetric:geominimal} and Proposition \ref{homogeneous:degree}, one has  
\begin{eqnarray*}
\frac{\tilde{G}_p(K)}{\tilde{G}_p(\E)}\geq \left(\frac{|\polar|}{|\E^\circ|}\right)^{\frac{p-n}{n+p}} \geq 1\ \ \  \Rightarrow\ \ \  {\tilde{G}_p(K)}\geq {\tilde{G}_p(\E)}.
\end{eqnarray*}  Clearly equality holds if $K=\E$. On the other hand, if $K\subsetneq \E$, then $\E^\circ\subsetneq \polar$, which implies $|\E^\circ |<|\polar|$ and hence ${\tilde{G}_p(K)}> {\tilde{G}_p(\E)}$. 

\vskip 2mm 
The following result compares the $L_p$ geominimal surface area with the $p$ surface area $S_p(K)=nV_p(K, \ball)$. Note that $S_p(\ball)=\tilde{G}_p(\ball)=n|\ball|=n\o_n$.   
\bp Let $K\in \cK_0$ be a convex body with the origin in its interior. \\
(i). For $p\geq 0$, one has $$\frac{\tilde{G}_p(K)}{\tilde{G}_p(\ball)}\leq \left(\frac{S_p(K)}{S_p(\ball)}\right)^{\frac{n}{n+p}}.$$
(ii). For $-n\neq p<0$, one has $$\frac{\tilde{G}_p(K)}{\tilde{G}_p(\ball)}\geq \left(\frac{S_p(K)}{S_p(\ball)}\right)^{\frac{n}{n+p}}.$$
\ep
\noindent {\bf Proof.} (i). For $p\geq 0$, by formula (\ref{Lp:geominimal:Surface:Area:positive}), one gets \begin{eqnarray*}
\tilde{G}_p(K)&=&\inf_{Q\in \cK_0} \left\{n\ V_p(K, Q) ^{\frac{n}{n+p}}\
  |Q^\circ|^{\frac{p}{n+p}}\right\}\\ &\leq& n\ [V_p(K, \ball)] ^{\frac{n}{n+p}}\
    |\ball|^{\frac{p}{n+p}}\\&=&(n\o_n)^{\frac{p}{n+p}} [S_p(K)]^{\frac{n}{n+p}}.
\end{eqnarray*} 
Dividing both sides by $n\o_n=\tilde{G}_p(\ball)=S_p(\ball)$, one gets the desired inequality.
\vskip 2mm \noindent 
(ii). For $-n\neq p<0$, by formula (\ref{Lp:geominimal:Surface:Area:negative}), one gets \begin{eqnarray*}
\tilde{G}_p(K)&=&\sup_{Q\in \cK_0} \left\{n\ V_p(K, Q) ^{\frac{n}{n+p}}\
  |Q^\circ|^{\frac{p}{n+p}}\right\}\\ &\geq& n\ [V_p(K, \ball)] ^{\frac{n}{n+p}}\
    |\ball|^{\frac{p}{n+p}}\\&=&(n\o_n)^{\frac{p}{n+p}} [S_p(K)]^{\frac{n}{n+p}}.
\end{eqnarray*} 
Dividing both sides by $n\o_n=\tilde{G}_p(\ball)=S_p(\ball)$, one gets the desired inequality.


\section{Cyclic inequalities and monotonicity of $\tilde{G}_p(\cdot)$} \label{section 5}

\bt\label{cyclic}
Let $K\in \cK_0$ be a convex body with the origin in its interior. 
\vskip 2mm \noindent
(i). If $-n<t<0<r<s$ or $-n<s<0<r<t$, then
\begin{equation}\label{i-2} \nonumber
\tilde{G}_r(K)\leq \tilde{G}_t(K) ^{\frac{(r-s)(n+t)}{(t-s)(n+r)}}
\tilde{G}_s (K)^{\frac{(t-r)(n+s)}{(t-s)(n+r)}}.
\end{equation}

\noindent (ii). If $-n<t<r<s<0$ or $-n<s<r<t<0$, then
\begin{equation}\label{i-3} \nonumber
\tilde{G}_r(K)\leq \tilde{G}_t(K)^{\frac{(r-s)(n+t)}{(t-s)(n+r)}}
\tilde{G}_s (K)^{\frac{(t-r)(n+s)}{(t-s)(n+r)}}.
\end{equation}

\noindent (iii). If $t<r<-n<s<0$ or $s<r<-n<t<0$, then
\begin{equation}\label{i-4} \nonumber
\tilde{G}_r(K)\geq \tilde{G}_t(K)^{\frac{(r-s)(n+t)}{(t-s)(n+r)}}
\tilde{G}_s (K)^{\frac{(t-r)(n+s)}{(t-s)(n+r)}}.
\end{equation}
\et
\vskip 2mm \noindent {\bf Proof.} Let $K, Q\in \cK_0$. We claim that, for all $r, s, t\in \bbR$  satisfying $0<\frac{t-r}{t-s}<1$,  
\be nV_r(K, Q)\leq [nV_s(K, Q)]^{\frac{t-r}{t-s}} \ [nV_t(K, Q)]^{\frac{r-s}{t-s}}.\label{mixed:p:surface:holder}\ee In fact, by H\"{o}lder's inequality (see \cite{HLP}) and $0<\frac{t-r}{t-s}<1$, \begin{eqnarray*}
nV_r(K, Q)&=&\int_{S^{n-1}} h_Q(u)^rh_K(u)^{1-r}\,dS(K, u)\\&=& \int_{S^{n-1}} [h_Q(u)^sh_K(u)^{1-s}]^{\frac{t-r}{t-s}} [h_Q(u)^th_K(u)^{1-t}]^{\frac{r-s}{t-s}}\,dS(K, u)\\ &\leq& 
\left(\int_{S^{n-1}} h_Q(u)^sh_K(u)^{1-s}\,dS(K, u)\right)^{\frac{t-r}{t-s}} \left(\int_{S^{n-1}} h_Q(u)^th_K(u)^{1-t}\,dS(K, u)\right)^{\frac{r-s}{t-s}}
\\&=& [nV_s(K, Q)]^{\frac{t-r}{t-s}} \ [nV_t(K, Q)]^{\frac{r-s}{t-s}}.
\end{eqnarray*} (i). Suppose that $-n<t<0<r<s$, which clearly implies $0<\frac{t-r}{t-s}<1$. First, we have \begin{eqnarray*}
\tilde{G}_t(K) ^{\frac{(r-s)(n+t)}{(t-s)(n+r)}}&=&\left\{ \sup_{Q_1\in \cK_0} [n V_t(K, Q_1)^{\frac{n}{n+t}}|Q_1^\circ|^{\frac{t}{n+t}}]\right\}^{\frac{(r-s)(n+t)}{(t-s)(n+r)}} \\&=&\sup_{Q_1\in \cK_0} \left\{[n V_t(K, Q_1)^{\frac{n}{n+t}}|Q_1^\circ|^{\frac{t}{n+t}}]^{\frac{(r-s)(n+t)}{(t-s)(n+r)}}\right\},
\end{eqnarray*} due to ${\frac{(r-s)(n+t)}{(t-s)(n+r)}}>0$. Similarly, as ${\frac{(t-r)(n+s)}{(t-s)(n+r)}}>0$, one has 
\begin{eqnarray*}\tilde{G}_s (K)^{\frac{(t-r)(n+s)}{(t-s)(n+r)}} &=& \inf_{Q\in \cK_0}\left\{[nV_s(K, Q)^{\frac{n}{n+s}}|Q^\circ|^{\frac{s}{n+s}}]^{\frac{(t-r)(n+s)}{(t-s)(n+r)}}\right\}.
\end{eqnarray*} By inequality (\ref{mixed:p:surface:holder}) and $\frac{n}{n+r}>0$, one has, $\forall Q\in \cK_0$,
\begin{eqnarray}
 \tilde{G}_r(K)&=&\inf_{Q_0\in \cK_0}\{ nV_r(K, Q_0) ^{\frac{n}{n+r}}\
  |Q_0^\circ|^{\frac{r}{n+r}}\}\leq nV_r(K, Q) ^{\frac{n}{n+r}}\
    |Q^\circ|^{\frac{r}{n+r}} \nonumber \\  
  &\leq& \big\{ [nV_s(K, Q)^{\frac{n}{n+s}}|Q^\circ|^{\frac{s}{n+s}}]^{\frac{(t-r)(n+s)}{(t-s)(n+r)}} \big\} \ \big\{[n V_t(K, Q)^{\frac{n}{n+t}}|Q^\circ|^{\frac{t}{n+t}}]^{\frac{(r-s)(n+t)}{(t-s)(n+r)}} \big\} \nonumber\\  &\leq&  \big\{[nV_s(K, Q)^{\frac{n}{n+s}}|Q^\circ|^{\frac{s}{n+s}}]^{\frac{(t-r)(n+s)}{(t-s)(n+r)}} \big\} \sup_{Q_1\in \cK_0} \left\{[n V_t(K, Q_1)^{\frac{n}{n+t}}|Q_1^\circ|^{\frac{t}{n+t}}]^{\frac{(r-s)(n+t)}{(t-s)(n+r)}}\right\}\nonumber\\ &=& \tilde{G}_t(K)  ^{\frac{(r-s)(n+t)}{(t-s)(n+r)}} \big\{
[nV_s(K, Q)^{\frac{n}{n+s}}|Q^\circ|^{\frac{s}{n+s}}]^{\frac{(t-r)(n+s)}{(t-s)(n+r)}}\big\}. \label{cyclic:-1-1}
\end{eqnarray}  
Taking the infimum over $Q\in \cK_0$ in inequality (\ref{cyclic:-1-1}), one gets
\begin{eqnarray*}
\tilde{G}_r(K)&\leq&  \tilde{G}_t(K)  ^{\frac{(r-s)(n+t)}{(t-s)(n+r)}} 
\inf_{Q\in \cK_0}\left\{[nV_s(K, Q)^{\frac{n}{n+s}}|Q^\circ|^{\frac{s}{n+s}}]^{\frac{(t-r)(n+s)}{(t-s)(n+r)}}\right\}\\&=&\tilde{G}_t(K) ^{\frac{(r-s)(n+t)}{(t-s)(n+r)}}
\tilde{G}_s (K)^{\frac{(t-r)(n+s)}{(t-s)(n+r)}}.
\end{eqnarray*}   The case $-n<s<0<r<t$ follows immediately by switching the roles of $t$ and $s$. 

\vskip 2mm \noindent (ii). Suppose that $-n<t<r<s<0$, which clearly implies $0<\frac{t-r}{t-s}<1$.  Then \begin{eqnarray*}  \tilde{G}_s (K) ^{\frac{(t-r)(n+s)}{(t-s)(n+r)}}&=&\left\{\sup_{Q\in \cK_0} [nV_s(K, Q)^{\frac{n}{n+s}}|Q^\circ|^{\frac{s}{n+s}}]\right\} ^{\frac{(t-r)(n+s)}{(t-s)(n+r)}} \\&=& \sup_{Q\in \cK_0} \left\{ [nV_s(K, Q)^{\frac{n}{n+s}}|Q^\circ|^{\frac{s}{n+s}}]^{\frac{(t-r)(n+s)}{(t-s)(n+r)}}\right\}, 
\end{eqnarray*} due to $\frac{(t-r)(n+s)}{(t-s)(n+r)}>0$.  Similarly, by  ${\frac{(r-s)(n+t)}{(t-s)(n+r)}}>0$, \begin{eqnarray*}
 \tilde{G}_t(K)  ^{\frac{(r-s)(n+t)}{(t-s)(n+r)}}&=& \sup_{Q\in \cK_0} \left\{ [n V_t(K, Q)^{\frac{n}{n+t}}|Q^\circ|^{\frac{t}{n+t}}]^{\frac{(r-s)(n+t)}{(t-s)(n+r)}} \right\}.
\end{eqnarray*} By inequality (\ref{mixed:p:surface:holder}) and $\frac{n}{n+r}>0$, one has, $\forall Q\in \cK_0$,  $$ nV_r(K, Q) ^{\frac{n}{n+r}}\
  |Q^\circ|^{\frac{r}{n+r}} \leq \big\{[nV_s(K, Q)^{\frac{n}{n+s}}|Q^\circ|^{\frac{s}{n+s}}]^{\frac{(t-r)(n+s)}{(t-s)(n+r)}} \big\}\big\{  [n V_t(K, Q)^{\frac{n}{n+t}}|Q^\circ|^{\frac{t}{n+t}}]^{\frac{(r-s)(n+t)}{(t-s)(n+r)}}\big\}.$$ Taking the supremum over $Q\in \cK_0$, one gets 
\begin{eqnarray*}
 \tilde{G}_r(K)\!\!\!&=&\!\!\!\sup_{Q\in \cK_0}\{ nV_r(K, Q) ^{\frac{n}{n+r}}\
  |Q^\circ|^{\frac{r}{n+r}}\}\\ \!\!\! &\leq&\!\!\!  \sup_{Q\in \cK_0} \left\{ [nV_s(K, Q)^{\frac{n}{n+s}}|Q^\circ|^{\frac{s}{n+s}}]^{\frac{(t-r)(n+s)}{(t-s)(n+r)}}\right\}  \sup_{Q\in \cK_0} \left\{ [n V_t(K, Q)^{\frac{n}{n+t}}|Q^\circ|^{\frac{t}{n+t}}]^{\frac{(r-s)(n+t)}{(t-s)(n+r)}} \right\}    \\\!\!\! &=&\!\!\! \tilde{G}_s (K) ^{\frac{(t-r)(n+s)}{(t-s)(n+r)}} \tilde{G}_t(K)  ^{\frac{(r-s)(n+t)}{(t-s)(n+r)}}.
\end{eqnarray*} The case $-n<s<r<t<0$ follows immediately by switching the roles of $t$ and $s$. 

\vskip 2mm \noindent (iii). Suppose that $t<r<-n<s<0$, which clearly implies $0<\frac{t-r}{t-s}<1$. Then \begin{eqnarray*} \tilde{G}_s(K)  ^{\frac{(t-r)(n+s)}{(t-s)(n+r)}}&=& \left\{\sup_{Q_1\in \cK_0} [nV_s(K, Q_1)^{\frac{n}{n+s}}|Q_1^\circ|^{\frac{s}{n+s}}]\right\}^{\frac{(t-r)(n+s)}{(t-s)(n+r)}}\\&=& \inf_{Q_1\in \cK_0} \left\{[nV_s(K, Q_1)^{\frac{n}{n+s}}|Q_1^\circ|^{\frac{s}{n+s}}]^{\frac{(t-r)(n+s)}{(t-s)(n+r)}}\right\}, 
\end{eqnarray*}  
due to $\frac{(t-r)(n+s)}{(t-s)(n+r)}<0$. Similarly, by ${\frac{(r-s)(n+t)}{(t-s)(n+r)}}>0$, one has, \begin{eqnarray*} \tilde{G}_t(K)  ^{\frac{(r-s)(n+t)}{(t-s)(n+r)}} &=&  \sup_{Q\in \cK_0}\left\{ [nV_t(K, Q)^{\frac{n}{n+t}}|Q^\circ|^{\frac{t}{n+t}}]^{\frac{(r-s)(n+t)}{(t-s)(n+r)}}\right\}.
\end{eqnarray*}  By inequality (\ref{mixed:p:surface:holder}) and $\frac{n}{n+r}<0$,  \begin{eqnarray*}
 \tilde{G}_r(K)&=&\sup_{Q_0\in \cK_0}\{ nV_r(K, Q_0) ^{\frac{n}{n+r}}\
  |Q_0^\circ|^{\frac{r}{n+r}}\} \geq  nV_r(K, Q) ^{\frac{n}{n+r}}\
    |Q^\circ|^{\frac{r}{n+r}} \\  
  &\geq&  \big\{[nV_s(K, Q)^{\frac{n}{n+s}}|Q^\circ|^{\frac{s}{n+s}}]^{\frac{(t-r)(n+s)}{(t-s)(n+r)}}\big\}  \  \big\{[n V_t(K, Q)^{\frac{n}{n+t}}|Q^\circ|^{\frac{t}{n+t}}]^{\frac{(r-s)(n+t)}{(t-s)(n+r)}}\big\}  \\  &\geq&  \big\{[n V_t(K, Q)^{\frac{n}{n+t}}|Q^\circ|^{\frac{t}{n+t}}]^{\frac{(r-s)(n+t)}{(t-s)(n+r)}} \big\}  \left\{\inf_{Q_1\in \cK_0} [nV_s(K, Q_1)^{\frac{n}{n+s}}|Q_1^\circ|^{\frac{s}{n+s}}]^{\frac{(t-r)(n+s)}{(t-s)(n+r)}}\right\}     \\&=& \tilde{G}_s(K)  ^{\frac{(t-r)(n+s)}{(t-s)(n+r)}} \big\{
[nV_t(K, Q)^{\frac{n}{n+t}}|Q^\circ|^{\frac{t}{n+t}}]^{\frac{(r-s)(n+t)}{(t-s)(n+r)}}\big\}. 
\end{eqnarray*} Therefore, taking the supremum over $Q\in \cK_0$, one gets
\begin{eqnarray*}
\tilde{G}_r(K)&\geq&  \tilde{G}_s(K)  ^{\frac{(t-r)(n+s)}{(t-s)(n+r)}} 
\sup_{Q\in \cK_0}\left\{[nV_t(K, Q)^{\frac{n}{n+t}}|Q^\circ|^{\frac{t}{n+t}}]^{\frac{(r-s)(n+t)}{(t-s)(n+r)}}\right\}\\&=&  \tilde{G}_s (K) ^{\frac{(t-r)(n+s)}{(t-s)(n+r)}}\tilde{G}_t(K)  ^{\frac{(r-s)(n+t)}{(t-s)(n+r)}}. 
\end{eqnarray*}   The case $s<r<-n<t<0$ follows immediately by switching the roles of $t$ and $s$.   

\vskip 2mm We now show the monotonicity of $\tilde{G}_p(\cdot)$.
\bt \label{monotone:geominimal:p--1} Let $K\in \cK_0$, and $p, q\neq 0$. \\ \noindent (i). If either $-n<q<p$ or $q<p<-n$, one has $$\bigg( \frac{\tilde{G}_q(K)}{n|K|}\bigg)^{\frac{n+q}{q}}\leq \bigg( \frac{\tilde{G}_p(K)}{n|K|}\bigg)^{\frac{n+p}{p}}.$$  \noindent (ii). If $q<-n<p$, one has $$\bigg( \frac{\tilde{G}_q(K)}{n|K|}\bigg)^{\frac{n+q}{q}}\geq \bigg( \frac{\tilde{G}_p(K)}{n|K|}\bigg)^{\frac{n+p}{p}}.$$ \et

\noindent {\bf Proof.} Recall that  $\tilde{G}_0(K)=n|K|$. Note that the statement of Theorem \ref{cyclic} does not include the cases $s=0$ or $r=0$ or $t=0$. However, from the proof of Theorem  \ref{cyclic}, one can easily see that cyclic inequalities still hold for (only) one of $r,s,t$ equal to $0$. 

\vskip 2mm \noindent (i). Note that $-n<q<p$ has three different cases: $0<q<p$, $-n<q<0<p$, and  $-n<q<p<0$.

\vskip 2mm \noindent Case 1: $0<q<p$. Put $t=0$, $r=q$ and $s=p$ in part (i) of Theorem \ref{cyclic},  then\begin{equation*}
\tilde{G}_q(K)\leq \tilde{G}_0(K) ^{\frac{(p-q)n}{p(n+q)}}
\tilde{G}_p (K)^{\frac{q(n+p)}{p(n+q)}}=\big(n|K| \big)^{\frac{(p-q)n}{p(n+q)}}
\tilde{G}_p (K)^{\frac{q(n+p)}{p(n+q)}}.
\end{equation*}
Dividing both sides by $n|K|$, one gets, as $q>0$, 
\begin{equation*}
\bigg(\frac{\tilde{G}_q(K)}{n|K|}\bigg)\leq 
\bigg(\frac{\tilde{G}_p (K)}{n|K|}\bigg)^{\frac{q(n+p)}{p(n+q)}}\Leftrightarrow \bigg( \frac{\tilde{G}_q(K)}{n|K|}\bigg)^{\frac{n+q}{q}}\leq \bigg( \frac{\tilde{G}_p(K)}{n|K|}\bigg)^{\frac{n+p}{p}}.
\end{equation*} 
Case 2: $-n<q<0<p$. Put $r=0$, $t=q$ and $s=p$ in part (i) of Theorem \ref{cyclic}, then
\begin{equation*}
\tilde{G}_0(K)=n|K| \leq  \tilde{G}_q(K)  ^{\frac{(n+q)p}{n(p-q)}}
 \tilde{G}_p (K) ^{\frac{q(n+p)}{(q-p)n}}.
\end{equation*}
Dividing both sides by $n|K|$, one gets, as $-n<q<0<p$, 
\begin{equation*}
\bigg(\frac{\tilde{G}_q(K)}{n|K|}\bigg)^{\frac{-p(n+q)}{(p-q)n}}\leq 
\bigg(\frac{\tilde{G}_p (K)}{n|K|}\bigg)^{\frac{q(n+p)}{(q-p)n}}\Leftrightarrow \bigg( \frac{\tilde{G}_q(K)}{n|K|}\bigg)^{\frac{n+q}{q}}\leq \bigg( \frac{\tilde{G}_p(K)}{n|K|}\bigg)^{\frac{n+p}{p}}.
\end{equation*} 
Case 3: $-n<q<p<0$. Put $s=0$, $t=q$ and $r=p$ in part (ii) of Theorem \ref{cyclic}, then  \begin{equation*}
 \tilde{G}_p(K)\leq  \tilde{G}_0(K)  ^{\frac{(q-p)n}{q(n+p)}}
  \tilde{G}_q (K) ^{\frac{p(n+q)}{q(n+p)}}=\big(n|K| \big)^{\frac{(q-p)n}{q(n+p)}}
   \tilde{G}_q (K) ^{\frac{p(n+q)}{q(n+p)}}.
 \end{equation*}
 Dividing both sides by $n|K|$, one gets, as $-n<p<0$, 
 \begin{equation*}
 \bigg(\frac{\tilde{G}_p(K)}{n|K|}\bigg)\leq 
 \bigg(\frac{\tilde{G}_q (K)}{n|K|}\bigg)^{\frac{p(n+q)}{q(n+p)}}\Leftrightarrow \bigg( \frac{\tilde{G}_q(K)}{n|K|}\bigg)^{\frac{n+q}{q}}\leq \bigg( \frac{\tilde{G}_p(K)}{n|K|}\bigg)^{\frac{n+p}{p}}.
 \end{equation*} 
 Case 4: $q<p<-n$. Put $s=0$, $t=q$ and $r=p$ in part (iii) of Theorem \ref{cyclic}, then
  \begin{equation*}
  \tilde{G}_p(K)\geq  \tilde{G}_0(K)  ^{\frac{(q-p)n}{q(n+p)}}
   \tilde{G}_q (K)^{\frac{p(n+q)}{q(n+p)}}=\big(n|K| \big)^{\frac{(q-p)n}{q(n+p)}}
   \tilde{G}_q (K)^{\frac{p(n+q)}{q(n+p)}}.
  \end{equation*}
  Dividing both sides by $n|K|$, one gets, as $p<-n$, 
  \begin{equation*}
  \bigg(\frac{\tilde{G}_p(K)}{n|K|}\bigg)\geq 
 \bigg(\frac{\tilde{G}_q (K)}{n|K|}\bigg)^{\frac{p(n+q)}{q(n+p)}}\Leftrightarrow \bigg( \frac{\tilde{G}_q(K)}{n|K|}\bigg)^{\frac{n+q}{q}}\leq \bigg( \frac{\tilde{G}_p(K)}{n|K|}\bigg)^{\frac{n+p}{p}}.
  \end{equation*} 
  (ii). Note that $q<-n<p$ has two different cases: $q<-n<0<p$ and $q<-n<p<0$. First, Proposition \ref{bounded by volume product}  and $\frac{q}{n+q}>0$ imply that $$ \frac{\tilde{G}_q(K)}{n|K|}\geq \bigg(\frac{|K^\circ|}{|K|}\bigg)^{\frac{q}{n+q}} \Leftrightarrow \bigg( \frac{\tilde{G}_q(K)}{n|K|}\bigg)^{\frac{n+q}{q}}\geq \frac{|K^\circ|}{|K|}.$$
  Similarly, for $p>0$, Proposition \ref{bounded by volume product}  and $\frac{p}{n+p}>0$ imply that $$ \frac{\tilde{G}_p(K)}{n|K|}\leq \left(\frac{|K^\circ|}{|K|}\right)^{\frac{p}{n+p}} \Rightarrow \bigg( \frac{\tilde{G}_p(K)}{n|K|}\bigg)^{\frac{n+p}{p}}\leq \frac{|K^\circ|}{|K|}\leq \bigg( \frac{\tilde{G}_q(K)}{n|K|}\bigg)^{\frac{n+q}{q}},$$ which concludes the case $q<-n<0<p$. For the case $q<-n<p<0$, by Proposition \ref{bounded by volume product}  and $\frac{p}{n+p}<0$, one has, $$ \frac{\tilde{G}_p(K)}{n|K|}\geq \left(\frac{|K^\circ|}{|K|}\right)^{\frac{p}{n+p}} \Rightarrow \bigg( \frac{\tilde{G}_p(K)}{n|K|}\bigg)^{\frac{n+p}{p}}\leq \frac{|K^\circ|}{|K|}\leq \bigg( \frac{\tilde{G}_q(K)}{n|K|}\bigg)^{\frac{n+q}{q}}.$$

 \noindent {\bf Remark.} In particular, if $p=1$, then for all $-n<q<0$ or $0<q<1$, \be \bigg( \frac{\tilde{G}_q(K)}{n|K|}\bigg)^{\frac{n+q}{q}}\leq \bigg( \frac{\tilde{G}_1(K)}{n|K|}\bigg)^{{n+1}}\label{comparison:1}.\ee 
 Now we can prove the following result which removes the centroid (or Santal\'{o} point) condition of $K$ in Theorem \ref{isoperimetric:geominimal}. 
 \bc \label{Lp:p:positive:symmetrization} Let $K\in\cK_0$ be a convex body with the origin in its interior. 
\\
\noindent (i). For $p\in (0, 1)$, the $L_p$
geominimal surface area attains the maximum at and only at origin-symmetric ellipsoids, among all convex bodies with fixed volume. More
precisely, \begin{eqnarray*} \frac{\tilde{G}_p(K)}{\tilde{G}_p(B^n_2)}\leq
\left(\frac{|K|}{|\ball |}\right)^{\frac{n-p}{n+p}}, 
\end{eqnarray*} with equality if and only if $K$ is an origin-symmetric ellipsoid. \\ \noindent 
(ii). For $p\in (-n, 0)$, the $L_p$
geominimal surface area attains the minimum at and only at origin-symmetric ellipsoids, among all convex bodies with fixed volume. More
precisely, \begin{eqnarray*} \frac{\tilde{G}_p(K)}{\tilde{G}_p(B^n_2)}\geq
\left(\frac{|K|}{|\ball|}\right)^{\frac{n-p}{n+p}}, 
\end{eqnarray*} with equality if and only if $K$ is an origin-symmetric ellipsoid. \ec

\vskip 2mm \noindent {\bf Proof.} Recall that $\tilde{G}_p(\ball)=n|\ball|$ for all $p\neq -n$ and the classical geominimal surface area $\tilde{G}_1(K)$ is translation invariant, namely, for all interior point $z_0$ of $K$, one has $$\tilde{G}_1(K)=\tilde{G}_1(K-z_0).$$ 

\vskip 2mm \noindent (i). Let $K\in \cK_0$ have centroid at $z_0\in \bbR^n$.  By inequality (\ref{comparison:1}), one has, as $p\in (0, 1)$, \begin{eqnarray*}
 \frac{\tilde{G}_p(K)}{n|K|} &\leq& \bigg( \frac{\tilde{G}_1(K)}{n|K|}\bigg)^{\frac{p(n+1)}{n+p}}=\bigg( \frac{\tilde{G}_1(K-z_0)}{\tilde{G}_1(\ball)}\bigg)^{\frac{p(n+1)}{n+p}}\bigg( \frac{|\ball|}{|K|}\bigg)^{\frac{p(n+1)}{n+p}}. 
\end{eqnarray*} Now using Theorem \ref{isoperimetric:geominimal} for $\tilde{G}_1(K-z_0)$ (as $K-z_0\in \cK_c$), one gets, 
\begin{eqnarray*}
 \frac{\tilde{G}_p(K)}{n|K|} &\leq& \bigg( \frac{\tilde{G}_1(K-z_0)}{\tilde{G}_1(\ball)} \bigg)^{\frac{p(n+1)}{n+p}} \bigg( \frac{|\ball|}{|K|} \bigg)^{\frac{p(n+1)}{n+p}}\leq  \bigg( \frac{|K-z_0|}{|\ball|} \bigg)^{\frac{p(n-1)}{n+p}} \bigg( \frac{|\ball|}{|K|} \bigg)^{\frac{p(n+1)}{n+p}}\\ &=& \bigg( \frac{|K|}{|\ball|} \bigg)^{\frac{p(n-1)}{n+p}} \bigg( \frac{|\ball|}{|K|} \bigg)^{\frac{p(n+1)}{n+p}}= \bigg( \frac{|K|}{|\ball|} \bigg)^{\frac{-2p}{n+p}}. 
\end{eqnarray*} Hence, one has 
\begin{eqnarray*}
  \frac{\tilde{G}_p(K)}{\tilde{G}_p(\ball)}  \cdot  \frac{|\ball|}{|K|}&\leq& \bigg( \frac{|K|}{|\ball|} \bigg)^{\frac{-2p}{n+p}}\ \ \ \Leftrightarrow\ \ \    \frac{\tilde{G}_p(K)}{\tilde{G}_p(\ball)}   \leq\left( \frac{|K|}{|\ball|}\right)^{\frac{n-p}{n+p}}. 
\end{eqnarray*} Clearly, equality holds if $K$ is an origin-symmetric ellipsoid. On the other hand, to have equality, one needs to have equality for the affine isoperimetric inequality related to $\tilde{G}_1(\cdot)$. Therefore, as explained in Theorem \ref{isoperimetric:geominimal}, $K-z_0 \in \cK_c$ must be an origin-symmetric ellipsoid. 
Proposition \ref{translation:invariance:fail:ball} implies that $z_0$ must be equal to $0$ and hence $K$ is an origin-symmetric ellipsoid. 

\vskip 2mm \noindent (ii). Let $K\in \cK_0$ have centroid at $z_0\in \bbR^n$.  By inequality (\ref{comparison:1}), one has, as $p\in (-n, 0)$, \begin{eqnarray*}
 \frac{\tilde{G}_p(K)}{n|K|}  &\geq&  \bigg( \frac{\tilde{G}_1(K)}{n|K|} \bigg)^{\frac{p(n+1)}{n+p}}= \bigg( \frac{\tilde{G}_1(K-z_0)}{\tilde{G}_1(\ball)} \bigg)^{\frac{p(n+1)}{n+p}} \bigg( \frac{|\ball|}{|K|} \bigg)^{\frac{p(n+1)}{n+p}}. 
\end{eqnarray*} Now using Theorem \ref{isoperimetric:geominimal} for $\tilde{G}_1(K-z_0)$ (as $K-z_0\in \cK_c$), one gets, as $p\in (-n, 0)$, 
\begin{eqnarray*}
  \frac{\tilde{G}_p(K)}{n|K|}  &\geq&  \bigg( \frac{\tilde{G}_1(K-z_0)}{\tilde{G}_1(\ball)} \bigg)^{\frac{p(n+1)}{n+p}}  \bigg( \frac{|\ball|}{|K|} \bigg)^{\frac{p(n+1)}{n+p}}\geq  \bigg( \frac{|K-z_0|}{|\ball|} \bigg)^{\frac{p(n-1)}{n+p}}  \bigg( \frac{|\ball|}{|K|} \bigg)^{\frac{p(n+1)}{n+p}}\\ &=& \bigg( \frac{|K|}{|\ball|} \bigg)^{\frac{p(n-1)}{n+p}} \bigg( \frac{|\ball|}{|K|} \bigg)^{\frac{p(n+1)}{n+p}}= \bigg( \frac{|K|}{|\ball|} \bigg)^{\frac{-2p}{n+p}}. 
\end{eqnarray*} Hence, one has 
\begin{eqnarray*}
 \frac{\tilde{G}_p(K)}{\tilde{G}_p(\ball)}\cdot  \frac{|\ball|}{|K|}&\geq&  \bigg( \frac{|K|}{|\ball|} \bigg)^{\frac{-2p}{n+p}}\ \ \ \Leftrightarrow \ \ \  \bigg( \frac{\tilde{G}_p(K)}{\tilde{G}_p(\ball)} \bigg) \geq \bigg( \frac{|K|}{|\ball|} \bigg)^{\frac{n-p}{n+p}}. 
\end{eqnarray*} Clearly, equality holds if $K$ is an origin-symmetric ellipsoid. On the other hand, to have equality, one needs to have equality for the affine isoperimetric inequality related to $\tilde{G}_1(\cdot)$. Therefore, as explained in Theorem \ref{isoperimetric:geominimal}, $K-z_0 \in \cK_c$ must be an origin-symmetric ellipsoid. 
Proposition \ref{translation:invariance:fail:ball} implies that $z_0$ must be equal to $0$ and hence $K$ is an origin-symmetric ellipsoid.

\vskip 2mm 
 \noindent {\bf Remark.}   Comparing the condition on $K$ in Corollary
\ref{Lp:p:positive:symmetrization} with those in Theorem \ref{isoperimetric:geominimal}, here 
one does not require the centroid (or the Santal\'{o} point) of $K$ to be at the origin. Analogous results for the $L_p$ affine surface area were first noticed in \cite{Zhang2007} by Zhang and further strengthened in \cite{Ye2013}. The case $p=1$
corresponds to the classical affine isoperimetric inequality related to the classical geominimal surface area \cite{Petty1974, Petty1985}.

\vskip 2mm \noindent {\bf Acknowledgments.} This paper is supported
 by a NSERC grant and a start-up grant from the Memorial University
 of Newfoundland. The author is grateful to the reviewer
 for many valuable comments.

\vskip 2mm \noindent Deping Ye, \ \ \ {\small \tt deping.ye@mun.ca}\\
{\small \em Department of Mathematics and Statistics\\
   Memorial University of Newfoundland\\
   St. John's, Newfoundland, Canada A1C 5S7 }

\end{document}